\documentclass[hidelinks,onefignum,onetabnum]{siamart251216}

\usepackage{todonotes}

\usepackage{lipsum}
\usepackage{amsfonts}
\usepackage{graphicx}
\usepackage{epstopdf}
\usepackage{algorithmic}
\ifpdf
  \DeclareGraphicsExtensions{.eps,.pdf,.png,.jpg}
\else
  \DeclareGraphicsExtensions{.eps}
\fi

\newsiamremark{remark}{Remark}
\newsiamremark{hypothesis}{Hypothesis}
\crefname{hypothesis}{Hypothesis}{Hypotheses}
\newsiamthm{claim}{Claim}
\newsiamremark{fact}{Fact}
\crefname{fact}{Fact}{Facts}

\headers{Adaptive Krylov methods for low-rank exponential integrators}{Rico Weigel, Tom-Christian Riemer, Martin Stoll}

\title{Adaptive Krylov methods for low-rank exponential integrators}

\author{Rico Weigel\thanks{Chair of Scientific Computing, Department of Mathematics, Technische Universit\"at Chemnitz, 09107 Chemnitz, Germany 
  (\email{rico.weigel@mathematik.tu-chemnitz.de}, \email{martin.stoll@mathematik.tu-chemnitz.de}, \email{tom.riemer@mathematik.tu-chemnitz.de}).}
\and Tom-Christian Riemer\footnotemark[1]
\and Martin Stoll\footnotemark[1]}

\usepackage{amsopn}

\usepackage{bbm}
\usepackage{nicefrac}
\usepackage{enumitem}
\usepackage{relsize}
\usepackage[noadjust]{cite}
\usepackage{graphicx}
\usepackage{subcaption}

\newcommand{\rkexp}{R{\smaller K2EXPINT}}
\newcommand{\rkexptt}{R{\smaller K2EXPINT\mbox{-}TT}}
\newcommand{\rkfit}{R{\smaller KFIT}}
\newcommand{\kiops}{K{\smaller IOPS}}
\newcommand{\kiopstt}{K{\smaller IOPS\mbox{-}TT}}
\newcommand{\tttext}{TT}
\newcommand{\ttbox}{TT-Toolbox}
\newcommand{\agmg}{A{\smaller GMG}}
\newcommand{\iom}{I{\smaller OM}}
\newcommand{\amen}{A{{\smaller ME}n}}
\newcommand{\ttmals}{TT\mbox{-}M{{\smaller ALS}}}
\newcommand{\ttgmres}{TT\mbox{-}G{{\smaller MRES}}}
\newcommand{\etderk}{E{\smaller TD3RK}}
\newcommand{\swt}{S{\smaller W2}}
\newcommand{\krogstad}{Krogstad4}
\newcommand{\iga}{IgA}
\newcommand{\fdm}{F{\smaller DM}}
\newcommand{\matlab}{M{\smaller ATLAB}}
\newcommand{\nurbs}{N{\smaller URBS}}
\newcommand{\svd}{S{\smaller VD}}
\newcommand{\hosvd}{H{\smaller OSVD}}

\newcommand{\C}{\mathbb{C}}
\newcommand{\R}{\mathbb{R}}
\newcommand{\N}{\mathbb{N}}
\newcommand{\one}{\mathbbm{1}}
\newcommand{\nset}{n_1,\ldots, n_d}

\newcommand{\iset}{i_1,\ldots,i_d}
\newcommand{\jset}{j_1,\ldots,j_d}
\newcommand{\ijset}{i_1, j_1, \ldots, i_d, j_d}
\newcommand{\jiset}{j_1, i_1, \ldots, j_d, i_d}

\newcommand{\tA}{\textbf{A}}

\newcommand{\tG}{\textbf{G}}

\newcommand{\tI}{\textbf{I}}

\newcommand{\tU}{\textbf{U}}

\newcommand{\tT}{\textbf{T}}
\newcommand{\tc}{\textbf{c}}
\newcommand{\tv}{\textbf{v}}
\newcommand{\tw}{\textbf{w}}
\newcommand{\tq}{\textbf{q}}
\newcommand{\tzero}{\textbf{0}}
\newcommand{\tvec}{\text{\normalfont{vec}\;\!}}
\newcommand{\tmat}{\text{\normalfont{mat}\;\!}}

\newcommand{\tF}{\textbf{F}}
\newcommand{\tg}{\textbf{g}}
\newcommand{\tu}{\textbf{u}}
\newcommand{\ta}{\textbf{a}}
\newcommand{\tb}{\textbf{b}}

\newcommand{\mm}{\textbf{m}}
\newcommand{\mn}{\textbf{n}}

\newcommand{\mmn}{\textbf{m}\times\textbf{n}}

\renewcommand{\mp}{\textbf{p}}
\newcommand{\nmodes}{n_1\times\ldots\times n_d}

\newcommand{\algorithmicinput}[1]{%
  \begin{description}[leftmargin=1.6cm, style=sameline]
    \item[\textbf{Input:}] #1
  \end{description}
}
\newcommand{\algorithmicoutput}[1]{%
  \begin{description}[leftmargin=1.6cm, style=sameline]
    \item[\textbf{Output:}] #1
  \end{description}
}

\newcommand{\hp}{\hphantom{}}

\ifpdf
\hypersetup{
  pdftitle={Adaptive Krylov methods for low-rank exponential integrators},
  pdfauthor={Rico Weigel, Martin Stoll and Tom-Christian Riemer}
}
\fi

\begin{document}

\maketitle

\begin{abstract}
    Differential equations arise in numerous applications, particularly within scientific and technical contexts. Systems of stiff, time-dependent ordinary differential equations constitute the focus of this work. Exponential integrators are designed to solve such equations by integrating the linear part exactly, while simultaneously approximating the nonlinear part through a linear combination of $\varphi$-functions. By utilizing an augmented stiffness matrix, state-of-the-art methods like \kiops\ and \rkexp\ solve the linear part and evaluate linear combinations of $\varphi$-functions for the nonlinear part in a single step, effectively reducing the computational effort to a single matrix exponential evaluation. However, these classical approaches assume that the system is represented by matrices and vectors, potentially not utilizing the underlying high-dimensional structure. Tensors address this limitation and offer significant storage efficiency through well-established decompositions like the Tensor Train (\tttext) format. This work provides a general framework for solving stiff, time-dependent systems directly within the \tttext\ format. Specifically, \kiopstt\ and \rkexptt\ are developed as extensions of the original \kiops\  and \rkexp\ algorithms. This involves reformulating the scheme of explicit exponential Runge-Kutta integrators for tensors and augmenting the stiffness tensor to compute linear combinations of $\varphi$-functions acting on tensors through a single evaluation of the exponential function using Krylov subspace methods. Furthermore, it is shown that the underlying theory of the matrix methods remains valid, thereby enabling the transfer of key theorems to the tensor case. Numerical experiments confirm significant speed-ups for \kiopstt\ and \rkexptt\ in low-rank scenarios compared to their classical counterparts.
\end{abstract}



\section{Introduction}\label{c1}
Differential equations arise in numerous applications and computing their solutions is often far from trivial. Traditional numerical methods typically represent systems of differential equations using matrix methods, whereby their solutions are expressed as vectors. Often these vector and matrix formulations do not entirely reflect the mathematical structure, especially in the higher-dimensional case. Furthermore, as all entries of a vector are typically stored explicitly, the computational effort increases exponentially for high-dimensional data. To preserve the mathematical structure while still reducing the storage requirements, low-rank tensor methods can be used, with well-established decompositions such as the Tensor Train (\tttext) format \cite{tt_oseledets, TT}.

In the context of using tensor methods for solving differential equations, methodological efforts center on two approaches that are also used in combination. The first approach attempts to integrate or update the individual components of a tensor decomposition, such as the \tttext\ cores in the case of the \tttext\ format, as separately as possible \cite{ceruti2020, ceruti2020.2, koch2010, kressner2010, lubich2015}. The second approach aims at addressing specific problems or problem structures \cite{kressner2010, cassini2023, munozmatute2022, croci2023, paeckel2019, lubich2015, li2019, vosidje2017}. For instance, it is frequently assumed that the linear operator is representable as a Kronecker sum.

In the case of stiff systems, the differential equation often consists of a stiff linear part, represented by the stiffness tensor matrix $\tA$, and a non-stiff semilinear part, represented by the tensor function $\tg$. Such problems arise, for example, in the discretization of semilinear parabolic differential equations on continuous domains \cite{RK2EXPINT}. Exponential integrators \cite{RK2EXPINT, expint, exprkint, KIOPS, PHIPM, etd3rk, krogstad4} have proven to be highly efficient solution methods in such cases, as they integrate the linear part $\tA$ exactly, while simultaneously approximating the nonlinear part. However, these integrators require the computation of linear combinations of $\varphi$-functions \cite{RK2EXPINT, expint, exprkint, KIOPS, PHIPM}, which can be resource-intensive. 

Regarding matrix methods for exponential integrators, Al-Mohy and Higham established a central theorem \cite{assembled_matrix} based on the earlier results of Saad \cite{saad} and Sidje \cite{sidje}. This theorem demonstrates that linear combinations of $\varphi$-functions can be evaluated by computing the matrix exponential of an augmented version of the stiffness matrix, reducing the computational cost to essentially a single evaluation of the matrix exponential. Leveraging these findings, the state-of-the-art methods \kiops\ \cite{KIOPS} and \rkexp\ \cite{RK2EXPINT} were developed to approximate the exponential of the augmented matrix via Krylov subspace methods and adaptively determine key parameters such as the Krylov space dimension.

To date, the application of exponential integrators to tensors in the \tttext\ format remains largely unexplored for systems lacking a more specific Kronecker sum structure. While Sidje’s work \cite{vosidje2017} seems to be the most related, it considers a linear differential equation, meaning that only the matrix exponential, and no additional $\varphi$-functions, needs to be evaluated. In contrast, this work provides a more general framework by deriving a method for solving systems of stiff, time-dependent differential equations, potentially featuring a nonlinear part, represented in the \tttext\ format. Specifically, the methods \kiopstt\ and \rkexptt\ are introduced as extensions of \kiops~\cite{KIOPS} and \rkexp\ \cite{RK2EXPINT}. In this approach, explicit exponential Runge-Kutta integrators~\cite{RK2EXPINT, expint, exprkint, KIOPS, PHIPM, etd3rk, krogstad4} are adapted for tensor-valued differential equations, involving the augmentation of the stiffness tensor to efficiently compute linear combinations of $\varphi$-functions. The resulting matrix exponential is then evaluated using Krylov subspace methods \cite{KIOPS, PHIPM, RK2EXPINT, saad_fab_and_error, lanczos, lanczos2, ruhe, ruhe2, ruhe3, ruhe4, guettel_diss, guettel_blockrat, faber}. Furthermore, it is shown that the underlying theory of the matrix-based versions remains applicable, thereby enabling the transfer of key theorems to the tensor case. A key advantage over tensor integrators that are limited to specific problem structures is that both the linear and nonlinear parts are not fixed in advance, allowing users to tailor them precisely to their specific needs.

However, Krylov methods for tensors have the disadvantage of requiring frequent multiplications of a tensor by a tensor matrix to find basis vectors, as well as numerous tensor sums for the orthogonalization process \cite{vosidje2017, amen}. This inherently increases the rank of the \tttext\ tensors, necessitating regular \tttext\ truncation (also referred to as rounding). Consequently, previous research on Krylov methods for \tttext\ tensors has been limited, with studies primarily focusing on solving linear systems of equations~\cite{tt_gmres, bucci2025}. The designs of \kiops\ and \rkexp, however, mitigate the frequency of required rounding in different ways. While \kiops\ uses an incomplete orthogonalization procedure for new basis vectors, thereby reducing the number of tensor sums, \rkexp\ aims to decrease the total number of iterations by utilizing a rational Krylov approach. These aspects motivate the investigation of Krylov methods for \tttext\ tensors despite the aforementioned drawbacks. In practice, numerical experiments confirm significant speedups of \kiopstt\ and \rkexptt\ in low-rank scenarios compared to the classical versions of \kiops\ and \rkexp.

The work is structured as follows: \Cref{c2} introduces basic concepts of tensors~\cite{TT} and the \tttext\ format \cite{tt_oseledets, TT}. \Cref{c3} presents classical exponential Runge-Kutta methods \cite{RK2EXPINT, expint, exprkint, KIOPS, PHIPM, etd3rk, krogstad4}, including the theorem for computing linear combinations of $\varphi$-functions \cite{RK2EXPINT, KIOPS, assembled_matrix}. In addition, this part presents the methods \kiops\ \cite{KIOPS} and \rkexp\ \cite{RK2EXPINT} as efficient approaches for evaluating the matrix exponential. All these theoretical foundations are incorporated into the development of \kiopstt\ and \rkexptt\ in \cref{c4}. Numerical experiments in \cref{c5} confirm the efficiency gains achieved through the use of the \tttext\ format. Finally, \cref{c6} draws major conclusions and provides an outlook on further development possibilities.

The complete source code and algorithmic implementations to reproduce the results of this work are publicly available at \url{https://github.com/riweig/ttexpint}.

\section{Tensors}\label{c2}
Tensors are used to store and process large amounts of data. This section provides the necessary background for defining the Tensor Train format and other essential notions used throughout this work. A comprehensive introduction to the theoretical framework of tensors can be found in the work of Gel\ss\ \cite{TT}. It also serves as the basis for this section and contains all concepts that are not explicitly defined in this paper but are used. Particularly in the case of ambiguities regarding the properties of the vectorization and matricization, consulting \cite{TT} is recommended.
\subsection{Basics}
A tensor can be seen as a $d$-dimensional array and will be denoted by a bold letter $\tT\in\C^{\nmodes}$. The numbers $\nset\in\N$ are called the modes. They are expressed in the mode set $\mn=(\nset)^T$. Therefore, it can be written $\tT\in\C^{\mn}$. Unless otherwise stated, it will be assumed that $\mm,\mn,\mp\in\N^d$ are $d$\mbox{-}dimensional mode sets. The element at position $\iset$ of a tensor is denoted by $\tT_{\iset}\in\C$ where $i_k\in\{1,\ldots,n_k\}$ for all $k=1,\ldots,d$. Although complex numbers are allowed as elements, this work mainly deals with real-valued tensors. The sum of two tensors and the product of a tensor with a scalar are defined element-wise, consistent with the standard definitions for matrices.

It will be necessary to apply linear operators $\tA:\C^\mn\rightarrow\C^\mm$, $\tT\mapsto\tA\tT$, on tensors. These can be expressed by a tensor $\tA\in\C^{\mmn}=\C^{m_1\times n_1\times\ldots\times m_d\times n_d}$. The action of $\tA$ on $\tT$ is defined as the product $\left(\tA\tT\right)\in\C^{\mm}$ where
\begin{equation*}
    \left(\tA\tT\right)_{\iset}:=\sum_{j_1=1}^{n_1}\cdots\sum_{j_d=1}^{n_d}\tA_{\ijset}\tT_{\jset}\text{.}
\end{equation*}
The linear map $\tA$ is called a tensor operator or, due to its action similar to the matrix-vector product, a tensor matrix. Each tensor $\tT\in\C^{\mn}$ has a tensor matrix representation $\widehat{\tT}\in\C^{\mn\times\one}$ where $\widehat{\tT}_{i_1, 1, \ldots, i_d, 1}:=\tT_{\iset}$. The transpose $\tA^T\in\C^{\mn\times\mm}$ of $\tA\in\C^{\mm\times\mn}$ is defined by $\tA^T_{\ijset}:=\tA_{\jiset}$. Accordingly, the adjoint $\tA^H$ is obtained by taking the transpose of $\tA$ and conjugating each entry. For a tensor $\tT\in\C^\mn$ it is $\tT^H:=\widehat{\tT}\vphantom{\tT}^H$. Using this, the $2$-norm of a tensor can be written as $\Vert\tT\Vert_2^2=\tT^H\tT$.

Let $\mm_{\text{prod}}$ and $\mn_{\text{prod}}$ denote the product of all modes in $\mm$ and $\mn$, respectively. The natural vectorization of a tensor $\tT\in\C^\mn$ reshapes it into a column vector $\tvec(\tT)\in\C^{\mn_{\text{prod}}}$ according to a bijective map $\psi_\mn$. The natural matricization transforms a tensor matrix $\tA\in\C^{\mm\times\mn}$ into a matrix $\tmat(\tA)\in\C^{\mm_{\text{prod}}\times\mn_{\text{prod}}}$ where the row index encapsulates the modes of $\mm$ and the column index spans the modes of $\mn$. The reindexing follows $\psi_\mn$ for the columns and an analogous bijection $\psi_\mm$ for the rows. As in \cite{TT}, these mappings are defined following the little-endian convention.

It is possible to generalize the concept of eigenvalues to the tensor framework. Define $\tzero_{\mn}\in\C^{\mn}$ as the zero tensor, i.e. the tensor containing only zeros. Let $\tA\in\C^{\mn\times\mn}$ be a tensor matrix. A scalar $\lambda\in\C$ that satisfies $\tA\tT=\lambda\tT$ for a tensor $\tT\in\C^{\mn}\setminus\{\tzero_\mn\}$ is called an eigenvalue of $\tA$ and $\tT$ is the corresponding eigentensor. It can be shown that every eigenvalue of $\tA$ is an eigenvalue of $\tmat(\tA)$ and vice versa.

Let $\tT\in\C^{m_1\times\ldots\times m_d}$ and $\tU\in\C^{n_1\times\ldots\times n_e}$ be two tensors. The tensor product ${\tT\otimes\tU\in\C^{m_1\times\ldots\times m_d\times n_1\times\ldots\times n_e}}$ is given by $\left(\tT\otimes\tU\right)_{i_1,\ldots,i_d,j_1,\ldots,j_e}:=\tT_{i_1,\ldots,i_d}\tU_{j_1,\ldots,j_e}$. This can be seen as a generalization of the outer product of two vectors. Furthermore, it can be shown that the tensor product of two tensor matrices represents the Kronecker product of their matricizations.

\subsection{Tensor Train format}

Consider a tensor in full format, i.e. every element is stored individually. The curse of dimensionality states that the size of the tensor, i.e. the number of its elements, grows exponentially with the number of dimensions $d$. This can be remedied by using tensor decompositions of low tensor rank. The Tensor Train format has gained attention because it mitigates the curse of dimensionality through a storage complexity that scales linearly with $d$ while ensuring algorithmic stability. The format originates from quantum physics, where it is known as the matrix product state representation, and was introduced to the mathematics community in~\cite{tt_oseledets}.
\begin{definition}
    A tensor $\tT\in\C^\mn$ or a tensor matrix $\tA\in\C^{\mm\times\mn}$ are represented in the Tensor Train (\tttext) format if
    \begin{align*}
        &\tT=\sum_{k_0=1}^{r_0}\cdots\sum_{k_d=1}^{r_d}\bigotimes_{i=1}^d\tT_{k_{i-1},:,k_i}^{(i)}
        &\text{or}
        &&\tA=\sum_{k_0=1}^{r_0}\cdots\sum_{k_d=1}^{r_d}\bigotimes_{i=1}^d\tA_{k_{i-1},:,:,k_i}^{(i)}
    \end{align*}
    where the numbers $r_0,\ldots,r_d\in\N$, $r_0=r_d=1$, are called the \tttext\ ranks and $\tT^{(i)}\in\C^{r_{i-1}\times n_i\times r_i}$ and $\tA^{(i)}\in\C^{r_{i-1}\times m_i\times n_i\times r_i}$ are the \tttext\ cores.
\end{definition}

Numerous operations can be efficiently executed within the \tttext\ format, including the addition of tensors, the action of tensor matrices on tensors, the computation of the $2$-norm and the solution of tensor-structured linear systems. However, the \tttext\ ranks typically increase during computations, necessitating regular \tttext\ rounding. To mitigate the resulting loss of accuracy, the rounding frequency must be carefully balanced against the expected growth of the \tttext\ ranks. Summation and multiplication represent two operations requiring careful adjustment. The ranks of a product of two \tttext\ matrices are bounded by the product of the corresponding ranks, i.e. they can increase drastically during multiplication. To maintain efficiency, the result of each product should be rounded directly. Regarding the sum of \tttext\ tensors, the problem is less critical, as the resulting ranks are bounded by the sum of the respective ranks of the summands. To avoid the loss of important information, the TTs are rounded every five sums or at the end of a sub-method in the implementation of the algorithms derived in this work.

For the implementation of low-rank algorithms, the \ttbox\ \cite{tt_toolbox} is used here. It contains various functions for working with \tttext\ tensors, such as initializing \tttext\ objects and performing calculations with them. Furthermore, it includes solvers for linear systems represented by \tttext\ tensors.

\section{Exponential integrators in the matrix case}\label{c3}
Consider the initial value problem 
\begin{equation}\label{eq_3_differential}
    \frac{\partial u(t)}{\partial t}=F(u(t))=-Au(t)+g(t,u(t)),\qquad u(0)=u_0
\end{equation}
from \cite{RK2EXPINT} where $u:[0,T]\to\R^n$ is an unknown function, $A\in\R^{n\times n}$ is a linear differential operator, $g:[0,T]\times\R^n\to\R^n$ is a semilinear function, i.e. it is generally nonlinear in $u$ but contains no derivatives of $u$, and $n\in\N$, $T\in\R_{\ge0}$. Furthermore, $u(0)=u_0\in\R^n$ is the initial condition.

A problem in solving this system is stiffness. While there is no clear definition of this term, a common property of stiff systems is that implicit numerical time integrators are superior to explicit ones \cite{RK2EXPINT, stiffness, lehrbuch_stiffness, lehrbuch_rungekutta, stiffness2}. This is an issue because implicit methods require the solution of (possibly nonlinear) systems of equations, while explicit methods require very small step sizes at high stiffness, and both strategies increase the computational cost \cite{lehrbuch_rungekutta}. Further possible characteristics of stiff problems of the form \eqref{eq_3_differential} are a large ratio of the absolute values of the largest and smallest eigenvalue or varying decay ratios of the solution's components \cite{RK2EXPINT, stiffness2}.

\subsection{Explicit exponential Runge-Kutta integrators}

To solve the differantial equation \cref{eq_3_differential}, the \textit{variation-of-constants} formula can be applied, yielding
\begin{equation*}
    u(t)=e^{-tA} u_0+\int_0^{t}e^{-(t-\tau)A} g(\tau, u(\tau))\text{ d}\tau
\end{equation*}
where $e^{-tA}$ denotes the matrix exponential of the matrix $-tA$ \cite{RK2EXPINT, expint, etd3rk, lehrbuch_rungekutta}. To enhance the quality of the solution, the time interval $[0,T]$ is divided into substeps $0=t_0<t_1<\ldots<t_k=T$. Furthermore, $h_i:=t_{i+1}-t_i$ for $i=0,\ldots,k-1$ and it is assumed that the solution $u_i:=u(t_i)$ is known for an $i\in\{0,\ldots,k-1\}$. When calculating the next step $u_{i+1}$, the last value $u_i$ can be considered as a new initial condition of \eqref{eq_3_differential}. As in \cite{RK2EXPINT, expint, etd3rk}, this leads to the expression
\begin{align}\label{eq_3_1_ui+1}
    u_{i+1}=u(t_i+h_i)=e^{-h_iA} u_i+\int_0^{h_i}e^{-(h_i-\tau)A} g(t_i+\tau, u(t_i+\tau))\text{ d}\tau\text{.}
\end{align}
Exponential integrators aim to integrate the linear part exactly, as it is often the source of stiffness, while approximating the nonlinear term using a linear combination of $\varphi$-functions \cite{RK2EXPINT, expint}.
\begin{definition}[\cite{RK2EXPINT, exprkint, PHIPM}]\label[definition]{def_3_1_phi}
    Let $k\in\N_0$ and $z\in\C$. The $k$-th $\varphi$-function is defined by the sum
    \begin{equation*}
        \varphi_{k}(z):=\sum_{j=0}^{\infty}\frac{z^j}{(j+k)!}\text{.}
    \end{equation*}
\end{definition}

The recurrence formula ${z\varphi_{k+1}(z)=\varphi_{k}(z)-\varphi_{k}(0)}$ can be derived by defining ${\varphi_0(z):=e^z}$ and $\varphi_{k}(0):=\nicefrac{1}{k!}$. It is possible to apply $\varphi$-functions to matrices ${A\in\C^{n\times n}}$ by the series definition or using $A\varphi_{k+1}(A)=\varphi_{k}(A)-\varphi_{k}(0)I_n=\varphi_{k+1}(A)A$ where ${I_n\in\C^{n\times n}}$ denotes the identity matrix. This can also be applied to tensor matrices ${\tA\in\C^{\mn\times\mn}}$ when using the identity tensor matrix $\tI_\mn\in\C^{\mn\times\mn}$ instead of $I_n$. 

Following the implementation of \rkexp~\cite{RK2EXPINT}, this work focuses on explicit exponential Runge-Kutta methods which belong to the class of exponential one-step methods \cite{expint}. They are designed to solve the initial value problem \cref{eq_3_differential} by 
dividing the integral of \eqref{eq_3_1_ui+1} into $s$ internal stages $0\le \gamma_1\le\ldots\le\gamma_s\le1$ and leveraging $\varphi$-functions to approximate it. The general explicit exponential Runge-Kutta scheme~\cite{RK2EXPINT, expint, lehrbuch_rungekutta} reads
\begin{equation}\label{eq_3_1_exprkint}
    \begin{aligned}
        u_{i+1}
        &=e^{-h_iA} u_i+h_i\sum_{j=1}^s\beta_j(-h_iA) G_{i,j}\text{,}\\
        U_{i,j}
        &=e^{-\gamma_jh_iA} u_i+h_i\sum_{k=1}^{j-1}\alpha_{j,k}(-h_iA) G_{i,k}\text{,}\\
        G_{i,k}
        &=g(t_i+\gamma_kh_i, U_{i,k})\text{.}
    \end{aligned}
\end{equation}
The coefficients $\beta_j(-h_iA)$ and $\alpha_{j,k}(-h_iA)$ represent $\varphi$-functions. As in \cite{RK2EXPINT}, these are chosen according to the two-stage integrator \swt\ \cite{sw2}, the three-stage method \etderk\ \cite{etd3rk} and the four-stage algorithm \krogstad\ \cite{krogstad4}.

At first glance, the evaluation of linear combinations of $\varphi$-functions does not appear to be very resource-efficient. In addition to the calculation of the matrix exponential for the $\varphi_0$-functions in each step, further $\varphi_k$ must be evaluated for $k=1,2,\ldots$, for example, by using the recursion formula. Fortunately, based on findings by Saad \cite{saad} and Sidje \cite{sidje}, it is shown by Al-Mohy and Higham \cite{assembled_matrix} that linear combinations of $\varphi$-functions can be computed efficiently by a slight enlargement of the matrix $A$ and a single evaluation of the matrix exponential \cite{RK2EXPINT}.
\begin{theorem}[\cite{RK2EXPINT, KIOPS, assembled_matrix}]\label[theorem]{thm_3_1_blockmatrix}
    Let $A\in\C^{n\times n}$, $C=[c_p,\ldots, c_1]\in\C^{n\times p}$ where $c_1,\ldots,c_p\in\C^n$ and $c_0\in\C^n$. Define $J_p\in\C^{p\times p}$ as a Jordan block to the eigenvalue 0. Furthermore, choose $h_i\in\R$ and set
    \begin{align*}
        &\Tilde{A}=\left[\begin{array}{cc}
        -A & C \\
        0 & J_p
        \end{array}\right]\in\C^{n+p\times n+p}
        &\text{and}
        &&\Tilde{c}=\left[\begin{array}{c}
             c_0 \\
             e_p 
        \end{array}\right]\in\C^{n+p}
    \end{align*}
    where $e_p=[0,\ldots,0,1]^T\in\C^{p}$. Then it is
    \begin{align*}
        &\left[e^{h_i\Tilde{A}}\right]_{1:n,n+p}=\sum_{k=1}^ph_i^k\varphi_k(-h_iA)c_k
        &\text{and}
        &&\left[e^{h_i\Tilde{A}}\Tilde{c}\right]_{1:n}=\sum_{k=0}^ph_i^k\varphi_k(-h_iA)c_k
    \end{align*}
\end{theorem}

A rearrangement of the terms in scheme \eqref{eq_3_1_exprkint} allows for taking full advantage of \cref{thm_3_1_blockmatrix}. A possible method is derived in \cite{expint}, which is also applied in the implementation of \cite{RK2EXPINT}. It remains to evaluate the matrix exponential efficiently, which is addressed in the next two subsections.

\subsection{\kiops}

A first approach of leveraging an adaptive Krylov method for computing the action of $\varphi$-functions on vectors to use them in exponential integrators is the algorithm 
{phipm} presented in \cite{PHIPM}. The \kiops\ algorithm is based on this idea and modifies it in three ways: it uses an augmented matrix instead of a time-stepping procedure to compute several actions of $\varphi$-functions at once, truncates the orthogonalization procedure in the Arnoldi algorithm and changes the adaptivity procedure \cite{KIOPS}. These modifications lead to an improvement in computational efficiency and stability as their results show. The contents of this subsection are based on the corresponding paper \cite{KIOPS}.

The primary focus in this subsection is the approximation of the term $e^{\Tilde{A}}\Tilde{c}$ where $\Tilde{A}\in\R^{n+p\times n+p}$ and $\Tilde{c}\in\R^{n+p}$. Krylov subspace methods represent an established strategy for this task.
\begin{definition}
    Let $\Tilde{A}\in\R^{n+p\times n+p}$ be a matrix, $\Tilde{c}\in\R^{n+p}$ be a vector and $m\in\N$. The (polynomial) Krylov subspace $\mathcal{K}_m(\Tilde{A},\Tilde{c})$ is defined as
    \begin{equation*}
        \mathcal{K}_m(\Tilde{A},\Tilde{c})=\text{\upshape span}\{\Tilde{c}, \Tilde{A}\Tilde{c}, \Tilde{A}^2\Tilde{c},\ldots, \Tilde{A}^{m-1}\Tilde{c}\}\text{.}
    \end{equation*}
\end{definition}

Ideally, one wants to find an orthonormal basis $V_m=[v_1,\ldots,v_m]\in\R^{n+p\times m}$ of $\mathcal{K}_m(\Tilde{A},\Tilde{c})$ and a matrix $H_m\in\R^{m\times m}$ representing $\Tilde{A}$ in the Krylov subspace, in short terms $\Tilde{A}V_m=V_mH_m$. But how to obtain $V_m$ and $H_m$? Since the assembled matrix $\Tilde{A}$ is non-symmetric by definition, the Lanczos iteration \cite{PHIPM, lanczos, lanczos2} cannot be applied. In this situation, a common approach is the Arnoldi iteration \cite{PHIPM, saad_fab_and_error}. For the inputs $\Tilde{A}$, $\Tilde{c}$, and the first basis vector $v_1=\nicefrac{\Tilde{c}}{\Vert \Tilde{c}\Vert_2}$, it produces an orthonormal basis $V_{m+1}=[v_1,\ldots,v_{m+1}]\in\R^{n+p\times m+1}$, a Hessenberg matrix $H_m\in\R^{m\times m}$ and an element $h_{m+1,m}\in\R$ satisfying the Arnoldi relation
\begin{align}\label{eq_3_2_arnoldi}
    \Tilde{A}V_m=V_mH_m+h_{m+1,m}v_{m+1}e_m^T\text{.}
\end{align}
An expensive step of the Arnoldi method is the orthogonalization of a new vector $v_{m+1}$ with respect to all previous vectors $v_1,\ldots,v_m$ in every iteration. To limit this, \kiops\ uses an incomplete orthogonalization procedure of length two, i.e. only the last two vectors $v_{m-1}$ and $v_m$ are considered during orthogonalization.

Based on the Arnoldi relation and as in \cite{KIOPS, PHIPM, saad_fab_and_error}, the approximation
\begin{align}\label{eq_3_2_approximation}
    e^{\Tilde{A}}v\approx\Vert \Tilde{c}\Vert_2V_me^{H_m}e_1\text{.}
\end{align}
can be derived. Assuming a small Krylov subspace size $m$, the cost of computing $e^{H_{m}}$ is negligible compared to the direct evaluation of $e^{A}$.

Two further aspects of the \kiops\ algorithm are not discussed here, as they are directly applicable to the tensor case without modification. These include, on the one hand, an a-posteriori error estimate and, on the other hand, an adaptivity method to determine the Krylov subspace size as well as an additional internal time step. Further details can be found in the original paper \cite{KIOPS}.

\subsection{\rkexp}

The method \rkexp\ \cite{RK2EXPINT} can be seen as a further developement of \kiops\ \cite{KIOPS}. It uses the adaptive procedure as a basis, but substitutes the polynomial Krylov algorithm by a  rational one. Rational Krylov methods are based on rational Krylov spaces, i.e. Krylov spaces whose elements are multiplied with the inverse of a polynomial, and were first introduced by Ruhe \cite{ruhe, ruhe2, ruhe3, ruhe4}. The substitution is justified by the finding that the exponential function can be approximated better by rational functions than by polynomials, see \cite{RK2EXPINT} for a brief discussion. This section covers the key aspects of \rkexp\ and is based on the corresponding paper \cite{RK2EXPINT}.

\begin{definition}[\cite{RK2EXPINT, guettel_diss}]\label[definition]{def_3_3_space_qm}
    Let $\Tilde{A}\in\C^{n+p\times n+p}$ be a matrix, $\Tilde{c}\in\C^{n+p}$ be a vector and $m\in\N$. Furthermore, define $\overline{\C}:=\C\cup\{\infty\}$. The rational Krylov subspace $\mathcal{Q}_m(\Tilde{A},\Tilde{c})$ is defined as
    \begin{equation*}
        \mathcal{Q}_m(\Tilde{A},\Tilde{c})=q_{m-1}(\Tilde{A})^{-1}\ \text{\upshape span}\{\Tilde{c}, \Tilde{A}\Tilde{c}, \Tilde{A}^2\Tilde{c},\ldots, \Tilde{A}^{m-1}\Tilde{c}\}=q_{m-1}(\Tilde{A})^{-1}\ \mathcal{K}_m(\Tilde{A},\Tilde{c})
    \end{equation*}
    where $q_{m-1}$ is the denominator polynomial given by $q_{m-1}(z)=\prod_{i=1,\xi_i\neq\infty}^{m-1}\left(z-\xi_i\right)$.
    The numbers $\xi_1,\ldots,\xi_{m-1}\in\overline{\C}$ are called the poles of $\mathcal{Q}_{m}(\Tilde{A},\Tilde{c})$ and must be different from the eigenvalues of $\Tilde{A}$.
\end{definition}

Once again, the aim is to represent $\Tilde{A}$ in a Krylov space, this time in a rational one. More precisely, the objective is to construct an orthonormal basis of $\mathcal{Q}_{m+1}(\Tilde{A},\Tilde{c})$, denoted by ${V_{m+1}=[v_1,\ldots,v_{m+1}]\in\C^{n+p\times m+1}}$, Hessenberg matrices $H_m,K_m\in\C^{m\times m}$ and $h_{m+1,m},k_{m+1,m}\in\C$ such that
\begin{align}\label{eq_3_3_ratarnoldi}
    \Tilde{A}V_mK_m+k_{m+1,m}\Tilde{A}v_{m+1}e_m^T=V_mH_m+h_{m+1,m}v_{m+1}e_m^T
\end{align}
where $\xi_mk_{m+1,m}=h_{m+1,m}$ if $\xi_m\neq\infty$ and $k_{m+1,m}=0$ else. Again, the first basis vector is assumed to be $v_1=\nicefrac{\Tilde{c}}{\Vert\Tilde{c}\Vert_2}$. Equation \cref{eq_3_3_ratarnoldi} is called the rational Arnoldi relation. The entries of $V_{m+1}, H_m, K_m$ and $h_{m+1,m},k_{m+1,m}$ can be computed by the rational Arnoldi method \cite{guettel_blockrat, guettel_toolbox}, which is based on Ruhe's algorithm \cite{ruhe, ruhe2, ruhe3, ruhe4}. This method may be regarded as an extension of the polynomial Arnoldi algorithm, where the next basis vector is obtained by solving a linear system rather than performing a matrix-vector multiplication. More precisely, if $\xi_i\neq\infty$, one has to solve
\begin{equation}\label{eq_3_3_lin_sys1}
    (\Tilde{A}-\xi_iI_{n+p})w=(\rho\Tilde{A}-\eta I_{n+p})V_iq
\end{equation}
where $\rho,\eta\in\{0,1\}$ are factors that are determined by a mobius transform of the pole $\xi_i$, and $q\in\C^i$ is the so-called continuation vector \cite{guettel_blockrat, ruhe4}. By multiplying both sides with $-1$ and defining $\Tilde{q}:=(\eta I_{n+p}-\rho\Tilde{A})V_iq$, the linear system can be written as
\begin{align*}
    (\xi_iI_{n+p}-\Tilde{A})w
    =\left[\begin{array}{cc}
        \xi_iI_n+A & -C \\
        0 & \xi_iI_p-J_p
    \end{array}\right]
    \left[\begin{array}{l}
        w_{1:n} \\
        \hp w_{n+1:n+p}
    \end{array}\right]
    =\left[\begin{array}{l}
        \Tilde{q}_{1:n} \\
        \hp\Tilde{q}_{n+1:n+p}
    \end{array}\right]\text{.}
\end{align*}
The lower block-row is assumed to contain only a few equations because $p$ should be small. Therefore, it can be solved directly. Based on $w_{n+1:n+p}$, the remaining components $w_{1:n}$ are computed using a suitable factorization or iterative method to solve the linear system
\begin{align}\label{eq_3_3_lin_sys2}
    (\xi_iI_n+A)w_{1:n}=\Tilde{q}_{1:n}+Cw_{n+1:n+p}\text{.}
\end{align}
For the numerical experiments in this work, the \agmg\ solver is used \cite{RK2EXPINT, agmg}, which is also considered in \cite{RK2EXPINT}. This choice is motivated by its iterative procedure, which is generally preferred for large-scale problems.

In the case of $\xi_i=\infty$ for all $i=1,\ldots,m$, the denominator polynomial vanishes in \cref{def_3_3_space_qm} and the rational Krylov space simplifies to a polynomial Krylov space. Accordingly, the rational Arnoldi algorithm reduces to the polynomial case, meaning that the solution of a linear system is bypassed. If at least $\xi_m=\infty$, \cref{eq_3_3_ratarnoldi} reduces to
\begin{equation*}
    \Tilde{A}V_mK_m=V_mH_m+h_{m+1,m}v_{m+1}e_m^T
\end{equation*}
and $K_m$ is invertible \cite{RK2EXPINT, guettel_fab, faber}. This leads to the approximation
\begin{align}\label{eq_3_3_approx}
    e^{h_i\Tilde{A}}\Tilde{c}
    \approx\Vert \Tilde{c}\Vert_2V_me^{h_iH_mK_m^{-1}}e_1\text{.}
\end{align}

Here, \rkexp \ uses the same adaptivity procedure as \kiops\ to choose a new Krylov space size. However, it is assumed that the maximum Krylov space size equals the number of poles provided by the user. When the poles are exhausted, \rkexp\ seems to want to continue with polynomial Krylov iterations, i.e. with infinity poles~\cite{RK2EXPINT, rk2expint_repo}. As part of this work, this method was revised so that \rkexp\ adjusts the step size like \kiops\ once the poles are exhausted. Another conceivable option would be a combination of both strategies, where infinity poles are added first and, when these are also exhausted, the step size is adjusted.

The a-posteriori error estimate of \rkexp\ remains unchanged within the scope of this work, and consequently requires no further discussion.

\section{Exponential integrators in the Tensor Train format}\label{c4}
Using tensors and tensor functions, the initial value problem \cref{eq_3_differential} can be written as
\begin{equation}\label{eq_4_differential}
    \frac{\partial \tu(t)}{\partial t}=\tF(\tu(t))=-\tA\tu(t)+\tg(t,\tu(t)),\qquad\tu(0)=\tu_0
\end{equation}
where $\tu:[0,T]\to\R^\mn$ is the unknown function, $\tA\in\R^{\mn\times\mn}$ is a linear differential operator represented by a tensor matrix, $\tg:[0,T]\times\R^\mn\to\R^\mn$ is a semilinear tensor function, and $\mn\in\N^d$, $T\in\R_{\ge0}$. Furthermore, $\tu(0)=\tu_0\in\R^\mn$ is the initial condition.

First, the definition of $\varphi$-functions for tensors is considered. The product of all modes in $\mn$ is again denoted by $\mn_{\text{prod}}$. By matricizing $\tA$ and using simple tensor manipulations (cf. \cite{TT}), the equation
\begin{equation*}
    \tmat\left(\varphi_{k}(\tA)\right)=\tmat\left(\sum_{j=0}^{\infty}\frac{1}{(j+k)!}\tA^j\right)=\sum_{j=0}^{\infty}\frac{1}{(j+k)!}\tmat(\tA)^j=\varphi_{k}(\tmat(\tA))
\end{equation*}
can be derived for $k\in\N_0$, and therefore it is also $\tmat(e^{\tA})=e^{\tmat(\tA)}$. Thus, the definitions of the tensor functions are consistent with their matrix counterpart. The ability to transform tensor matrices into matrices and tensors into vectors, and then to use matrix methods to solve tensor equations, allows explicit exponential Runge-Kutta integrators to be applied to tensor functions. By taking $\tA$, $\tu$ and $\tg$ as in \cref{eq_4_differential} and using the same time-stepping procedure as in \cref{c4}, scheme \eqref{eq_3_1_exprkint} can be reformulated for tensors as
\begin{equation}\label{eq_4_1_exprkint}
    \begin{alignedat}{2}
        \tu_{i+1}
        &=e^{-h_i\tA} \tu_i+h_i\sum_{j=1}^s\tb_j(-h_i\tA) \tG_{i,j}&\in\R^{\mn}\text{,} \\
        \tU_{i,j}
        &=e^{-\gamma_jh_i\tA} \tu_i+h_i\sum_{k=1}^{j-1}\ta_{j,k}(-h_i\tA) \tG_{i,k}\ &\in\R^{\mn}\text{,} \\
        \tG_{i,k}
        &=\tg(t_i+\gamma_kh_i, \tU_{i,k})&\in\R^{\mn}\text{.}
    \end{alignedat}
\end{equation}

It can be noticed that the $\varphi$-functions $\ta_{j,k}$ and $\tb_j$ now generate tensor matrices instead of matrices. To compute their linear combinations efficiently, \cref{thm_3_1_blockmatrix} should be applied again. In the matrix setting, $A\in\C^{n\times n}$, $C=[c_p,\ldots, c_1]\in\C^{n\times p}$ containing the vectors which are multiplied with the $\varphi$-functions, and $J_p\in\C^{p\times p}$ as a Jordan block to the eigenvalue 0 were combined to form
\begin{equation*}
    \Tilde{A}:=\left[\begin{array}{cc}
        -A & C \\
        0 & J_p
        \end{array}\right]\text{.}
\end{equation*}
But how to enlarge $\tA\in\C^{\mn\times\mn}$ as a tensor matrix? The first idea is to somehow directly enlarge the tensor matrix. There are two possible ways: Increasing one mode of $\mn$ or introducing a new dimension, i.e. adding a new element to $\mn$. Both ideas would lead to a substantial increase in tensor size, not to mention the problem of how to stuff the tensors $\tc_1,\ldots,\tc_p\in\C^\mn$ and the Jordan block $J_p$ into the assembled tensor matrix. It is even challenging to string the tensors $\tc_1,\ldots,\tc_p$ together. A possible solution is to consider the matricization, i.e. to use the matrix
\begin{equation}\label{eq_4_assembled_matrix}
    \Tilde{A}:=\left[\begin{array}{ccccc}
        -\tmat(\tA) & \vline & \tvec(\tc_p) & \cdots & \tvec(\tc_1) \\
        \hline
        0 & \vline & & J_p
        \end{array}\right]\text{.}
\end{equation}
This has the advantage that \cref{thm_3_1_blockmatrix} remains directly applicable, yielding
\begin{equation}\label{eq_4_assembled_exp}
\begin{aligned}
        \left[e^{h_i\Tilde{A}}\left[\begin{array}{c}
            \tvec(\tc_0) \\
            e_p
        \end{array}\right]\right]_{1:\mn_{\text{prod}}}
        =\tvec\left(\sum_{k=0}^ph_i^k\varphi_k(-h_i\tA)\tc_k\right)
    \end{aligned}
\end{equation}

where $\tc_0\in\C^\mn$ and leveraging properties of the vectorization (cf. \cite{TT}). As the vectorization reshapes the tensor according to a bijective map, it can be reversed and thus, it is possible to obtain the linear combination of $\varphi$-functions directly as a tensor. But it should be kept in mind that, since $\tA$ and $\tc_0,\ldots,\tc_p$ are assumed to be stored in the \tttext\ format, the computation of the matricization and vectorization must be avoided at all costs. Otherwise, all entries of the objects would be stored explicitly, which eliminates any efficiency advantages of the TT format. The next two subsections explore ways to treat the blocks of $\Tilde{A}$ individually, instead of assembling the matrix. They also examine methods to obtain the term $\sum_{k=0}^ph_i^k\varphi_k(-h_i\tA)\tc_k$ directly in tensor form, rather than as a vector.

\subsection{Modifications to \kiops}\label[subsection]{s_4_1}

The polynomial Krylov space $\mathcal{K}_m(\Tilde{A},\Tilde{c})$ and the incomplete orthogonalization method (\iom) of \kiops\ \cite{KIOPS} are considered here again. The goal is to apply the \iom\ directly to tensors and tensor matrices instead of reshaping them. In this context, the property of \kiops\ that $\Tilde{A}$ is not explicitly assembled, but its blocks are treated separately, can be leveraged. Let $\Tilde{A}$ be defined as in \cref{eq_4_assembled_matrix} and $\Tilde{c}:=\left[\tvec(\tc_0)^T, e_p^T\right]^T$ for $\tc_0\in\R^{\mn}$. 
In contrast to the \iom, the input parameters are modified as follows to obtain \cref{alg_4_2_iom}:

\begin{itemize}[label=$\circ$]
    \item The linear operator $\tA\in\R^{\mn\times\mn}$ remains in its original form as a tensor matrix instead of being matricized.
    \item The tensors $\tc_1,\ldots,\tc_p\in\R^{\mn}$ are no longer stacked column-wise and remain in their tensor form, too.
    \item Each former vector $v_i\in\R^{\mn_{\text{prod}}+p}$, where $i=1,\ldots,m+1$, and $\mn_{\text{prod}}$ is the product of all modes in the mode set $\mn$, is split into two parts: a tensor $\tv_i\in\R^\mn$ containing the first $\mn_{\text{prod}}$ entries and a vector $v_i^{(p)}\in\R^p$ containing the last $p$ entries. The vectors $v_i^{(p)}$ are assembled in the matrix ${V_{m+1}^{(p)}:=[v_1^{(p)},\ldots,v_{m+1}^{(p)}]\in\R^{p\times m+1}}$.
    \item It is $\tv_1=\nicefrac{\tc_0}{c_{\text{norm}}}$ and $v_1^{(p)}=\nicefrac{e_p}{c_{\text{norm}}}$ where $c_\text{norm}:=(\Vert\tc_0\Vert_2^2+1)^{\nicefrac{1}{2}}$ and $\tc_0\in\R^n$.
\end{itemize}

The output of \cref{alg_4_2_iom} satisfies the relation
\begin{align}
    -\tA\tv_i+\sum_{k=1}^p\tc_{p-k+1}\left[v_i^{(p)}\right]_k
    &=\sum_{k=1}^m\tv_k[H_m]_{k,i}+\delta_{i,m}h_{m+1,m}\tv_{m+1}\label{eq_4_2_output_iom1} \\
    J_pv_i
    &=V_m[H_m]_{:,i}+\delta_{i,m}h_{m+1,m}v^{(p)}_{m+1}\label{eq_4_2_output_iom2}
\end{align}
for $i=1,\ldots,m$ and $\delta_{i,m}$ as the Kronecker delta. Vectorizing both sides of formula~\cref{eq_4_2_output_iom1} yields
\begin{equation*}
    \text{\small$
    -\tmat(\tA)\tvec(\tv_i)+\sum_{k=1}^p\tvec(\tc_{p-k+1})\left[v_i^{(p)}\right]_k \\
    =\sum_{k=1}^m\tvec(\tv_k)[H_m]_{k,i}+\delta_{i,m}h_{m+1,m}\tvec(\tv_{m+1})\text{.}$}
\end{equation*}
Let $\Tilde{A}$ be defined as in \eqref{eq_4_assembled_matrix}, $v_i^{(\mn)}:=\tvec(\tv_i)$ for $i=1,\ldots,m+1$, and $V_m^{(\mn)}:=[v_1^{(\mn)},\ldots,v_m^{(\mn)}]$. Taking \eqref{eq_4_2_output_iom2} into account, the expression
\begin{align}\label{eq_4_2_tt_arnoldi}
    \Tilde{A}\left[\begin{array}{ccc}
        V_m^{(\mn)} \\
        V_m^{(p)}
    \end{array}\right]
    &=\left[\begin{array}{c} 
        V_m^{(\mn)} \\
        V_m^{(p)}
    \end{array}\right]H_m+h_{m+1,m}\left[\begin{array}{c} 
        v_{m+1}^{(\mn)} \\
        v_{m+1}^{(p)}
    \end{array}\right]e_m^T
\end{align}
can be derived, which is equivalent to the Arnoldi relation \eqref{eq_3_2_arnoldi}. The tensor version of the \iom\ is summarized in \cref{alg_4_2_iom}.

\begin{algorithm}[H]\label[algorithm]{alg_4_2_iom}
    \caption{Tensor version of the \iom}
    \algorithmicinput{$\tA\in\R^{\mn\times \mn}$, $\tc_1,\ldots,\tc_p\in\R^{\mn}$ and $\tv_1,\ldots,\tv_{m+1}\in\R^{\mn}$, ${v_1^{(p)},\ldots,v_{m+1}^{(p)}\in\R^{p}}$, $H\in\R^{m+1\times m}$ storing a Hessenberg matrix $H_m\in\R^{m\times m}$ and $h_{m+1,m}$ satisfying the Arnoldi relation \eqref{eq_4_2_tt_arnoldi} and $m_{\text{max}}$}
    \algorithmicoutput{$\tv_1,\ldots,\tv_{m_{\text{end}}}\in\R^{\mn}$, $v_1^{(p)},\ldots,v_{m_{\text{end}}}^{(p)}\in\R^{p}$ and, if happy\_breakdown is False, ${H\in\R^{m_\text{end}\times m_\text{end}-1}}$ where $m_{\text{end}}=m_{\text{max}}+1$ and else $H\in\R^{m_\text{end}\times m_\text{end}}$ where $m+1\le m_{\text{end}}\le m_{\text{max}}$}
    \begin{algorithmic}[1]
        \STATE{happy\_breakdown = False}
        \FOR{$i=m+1,\ldots,m_{\text{max}}$}
            \STATE{$\tw = -\tA\tv_i + \sum_{k=1}^p\tc_{p-k+1}v_i^{(p)}(k)$}
            \STATE{$w = J_pv_i^{(p)}$}
            \FOR{$j=\max(1,i-1),\ldots,i$}
                \STATE{$H(i,j) = \tv_i^T\tw+\left(v_i^{(p)}\right)^Tw$}
                \STATE{$\tw = \tw-H(i,j)\tv_i$}
                \STATE{$w = w-H(i,j)v_i^{(p)}$}
            \ENDFOR
            \STATE{$h_{\text{new}}=\Vert \tw\Vert_2^2+\Vert w\Vert_2^2$}
            \IF{$\sqrt{h_{\text{new}}} \approx 0$}
                \STATE{happy\_breakdown = True}
                \STATE{\textbf{break}}
            \ENDIF
            \STATE{$H(i+1,i)=\sqrt{h_{\text{new}}}$}
            \STATE{$\tv_{i+1}={\tw}/{H(i+1,i)}$}
            \STATE{$v_{i+1}^{(p)}={w}/{H(i+1,i)}$}
        \ENDFOR
    \end{algorithmic}
\end{algorithm}

With the above derivations in mind, approximation \cref{eq_3_2_approximation} can be used to obtain
\begin{equation*}
    e^{h_i\Tilde{A}}\left[\begin{array}{c} 
        c_0 \\
        e_p
    \end{array}\right]\approx\left(\Vert\tc_0\Vert_2^2+1\right)^{\nicefrac{1}{2}}\left[\begin{array}{c} 
        V_m^{(\mn)} \\
        V_m^{(p)}
    \end{array}\right]e^{h_iH_m}e_1
\end{equation*}
where $c_0:=\tvec(\tc_0)$ for $\tc_0\in\R^\mn$, and $\Tilde{A}$ is the matricized tensor matrix defined in~\cref{eq_4_assembled_matrix}. Considering all but the last $p$ rows and using formula \cref{eq_4_assembled_exp} yields
\begin{align*}
    \tvec\left(\sum_{k=0}^ph_i^k\varphi_k(-h_i\tA)\tc_k\right)
    \approx\tvec\left(c_\text{norm}\left(\sum_{i=1}^m\tv_i\left[e^{h_iH_m}\right]_{i,1}\right)\right)\text{.}
\end{align*}
Since the vectorization reshapes all tensors of this equation following the same bijective map, the approximation is finally obtained in tensor form as
\begin{align*}
    \sum_{k=0}^ph_i^k\varphi_k(-h_i\tA)\tc_k
    &\approx c_\text{norm}\left(\sum_{i=1}^m\tv_i\left[e^{h_iH_m}\right]_{i,1}\right)\text{.}
\end{align*}
In summary, it can be stated that, in the polynomial case, all calculations can be performed directly with tensors. The need for explicit matricizations or vectorizations is entirely avoided within this framework. Therefore, the algorithms can fully benefit from the use of the \tttext\ format. The method derived here is referred to as \kiopstt.

\subsection{Modifications to \rkexp}\label[subsection]{s_4_2}

This section focuses on rational Krylov spaces $\mathcal{Q}_m(\Tilde{A},\Tilde{c})$ and the rational Arnoldi iteration for \tttext\ tensors. The modification to the input of the rational Arnoldi algorithm are the same as in \Cref{s_4_1}. In contrast to \rkexp's matrix-based rational Krylov iteration \cite{RK2EXPINT}, the matrix $\Tilde{A}$ should not be assembled. This also means that $\tA\in\R^{\mn\times\mn}$ and $\tc_1,\ldots,\tc_p\in\R^\mn$ should not be matricized or vectorized. In particular, care must be taken when solving the equivalent to system \cref{eq_3_3_lin_sys1}. Consider its right-hand side 
\begin{equation*}
    \Tilde{q}:=(\eta I_{\mn_{\text{prod}}+p}-\rho\Tilde{A})V_iq\text{.}
\end{equation*} 
The last $p$ rows can be written as $\Tilde{q}^{(p)}:=(\eta I_{p}-\rho J_p)V_i^{(p)}q$ and therefore, the computation of $\Tilde{q}^{(p)}$ remains the same since all of its objects are standard matrices or vectors by definition. The first $\mn_{\text{prod}}$ rows of $\Tilde{q}$ are
\begin{align*}
    \Tilde{q}^{(\mn)}
    &:=\left[\begin{array}{cc}
        \eta I_{\mn_{\text{prod}}}+\rho A & -\rho C
    \end{array}\right]V_iq
    =(\eta I_{\mn_{\text{prod}}}+\rho A)V_i^{(\mn)}q-\rho CV_i^{(p)}q\text{.}
\end{align*}
Using that $q\in\R^i$ is a vector, $\Tilde{q}^{(\mn)}$ can be rewritten as
\begin{align*}
    \Tilde{q}^{(\mn)}
    &=\tvec\left((\eta\tI_\mn+\rho\tA)\sum_{j=1}^i\tv_jq_j-\rho\sum_{k=1}^p\tc_{p-k+1}\left[V_i^{(p)}q\right]_{k}\right)
    =:\tvec(\Tilde{\tq})\text{.}
\end{align*}

With this in mind, the focus can now be shifted to the left-hand side of the linear system. Define $w^{(p)}:=w_{\mn_{\text{prod}}+1:\mn_{\text{prod}}+p}$ and $w^{(\mn)}:=w_{1:\mn_{\text{prod}}}=:\tvec(\tw)$ where $\tw\in\R^\mn$. The last $p$ rows form again a small system
\begin{equation*}
    \xi_iI_p-J_pw^{(p)}=\Tilde{q}_i^{(p)},
\end{equation*}
which can be solved directly. There are no changes necessary because all objects are still vectors or matrices. The first $\mn_{\text{prod}}$ rows can be expressed by
\begin{align*}
    (\xi_iI_{\mn_{\text{prod}}}+A)w^{(\mn)}-Cw^{(p)}
    &=\tvec\left((\xi_i\tI_\mn+\tA)\tw\right)-\tvec\left(\sum_{k=1}^p\tc_{p-k+1}w^{(p)}_k\right)\text{.}
\end{align*}
Together with the right-hand side, this forms a linear system in the tensor framework:
\begin{align}\label{eq_4_3_tensor_system}
    (\xi_i\tI_\mn+\tA)\tw=\Tilde{\tq}+\sum_{k=1}^p\tc_{p-k+1}w^{(p)}_k\text{.}
\end{align}
By exploiting the tensor format of the input data, the solution $\tw$ can be computed efficiently. Since it is infeasible to calculate an exact solution for large tensors, the result should be approximated. Numerous solution methods exist for this task that utilize the \tttext\ format, e.g. the \ttmals\ solver \cite{mals}, the \ttgmres\ solver \cite{tt_gmres} or the \amen\ solver \cite{amen}. For the implementation, the \amen\ solver \cite{amen, tt_toolbox} is chosen as it represents one of the fastest methods \cite{amen, amen_comp} and delivers reliable results also for non-symmetric systems \cite{amen}.

The poles $\xi_i$, $i=1,\ldots,m$, of the rational Krylov space $\mathcal{Q}_{m+1}(\Tilde{A},\Tilde{c})$ must be examined too. It is assumed that these are not located within the spectrum of $\Tilde{A}$. Since $\Tilde{A}$ is an upper block triangular matrix, its spectrum is equal to the union of the spectra of the matrices $A$ and $J_p$. Since $J_p$ is consistently treated as a matrix in this section, its only eigenvalue is zero. For $A=\tmat(\tA)$, the eigenvalues of $A$ are equal to the eigenvalues of $\tA$, and therefore no complications are encountered here either.

When considering complex poles, the linear system \eqref{eq_4_3_tensor_system} becomes complex-valued. Numerical experiments reveal that the \amen\ solver encounters stability issues when solving such systems directly. To address this, a second variant of the rational Arnoldi algorithm for TTs is implemented. This alternative approach separates the real and imaginary parts of all involved TT tensors, thereby restricting all computations to the real domain. Consequently, the complex system \eqref{eq_4_3_tensor_system} is transformed into the real-valued block-tensor system
\begin{equation*}
    \left[\begin{array}{cc}
        \text{re}(\xi_i\tI_{\mn}+\tA) & -\text{im}(\xi_i\tI_{\mn}+\tA) \\
        \text{im}(\xi_i\tI_{\mn}+\tA) & \text{re}(\xi_i\tI_{\mn}+\tA)
    \end{array}\right]\left[\begin{array}{c}
        \text{re}(\tw) \\
        \text{im}(\tw) 
    \end{array}\right]=\left[\begin{array}{c}
        \text{re}(\Tilde{\tq}+\sum_{k=1}^p\tc_{p-k+1}w_k^{(p)}) \\
        \text{im}(\Tilde{\tq}+\sum_{k=1}^p\tc_{p-k+1}w_k^{(p)})
    \end{array}\right]
\end{equation*}
where $\text{re}(\cdot)$ and $\text{im}(\cdot)$ denote the real and imaginary parts of an object, respectively~\cite{cplx_lin_sys}. Generally, splitting an object into its real and imaginary components doubles its total number of entries. For tensors and tensor matrices, the splitting can be achieved by appending a new mode of size two to the existing mode set. In the TT format, this can be implemented at low additional cost by adding a new core. For instance, the transformed inputs are defined as ${\tA^{(\text{split})}:=\tA\otimes I_2}$ with $I_2\in\R^{2\times2}$ as the identity matrix, ${\tI_{\mn}^{(\text{split})}:=\tI_\mn\otimes[[\text{re}(\xi_i), \text{im}(\xi_i)]^T, [-\text{im}(\xi_i), \text{re}(\xi_i)]^T]}$, ${\tc_i^{(\text{split})}:=\tc_i\otimes[1, 0]^T}$ for $i=1,\ldots,p$ and $\tv_1^{(\text{split})}:=\tv_1\otimes[1, 0]^T$. In the final algorithm, this splitting naturally extends to the tensors $\Tilde{\tq}$, $\tw$ and $\tv_i$, $i=2,\ldots,m+1$.

Please keep in mind that in case of complex-valued poles, the computation of the product of a complex number with a TT tensor has to take care of the split of the real and complex parts of the TT tensors. Since the overall solution is assumed to be real valued, it can be extracted by fixing the last, additional dimension to one.

All other modifications to the rational Arnoldi iteration are quite similar to those made in \Cref{s_4_1} and will be assumed to hold from the previous discussion. The output of the modified procedure satisfies the rational Arnoldi relation
\begin{equation*}
    \begin{aligned}
        \Tilde{A}\left(\left[\begin{array}{ccc}
            V_m^{(\mn)} \\
            V_m^{(p)}
        \end{array}\right]K_m+k_{m+1,m}\left[\begin{array}{ccc}
            v_{m+1}^{(\mn)} \\
            v_{m+1}^{(p)}
        \end{array}\right]e_m^T\right)
        =\left[\begin{array}{c} 
            V_m^{(\mn)} \\
            V_m^{(p)}
        \end{array}\right]H_m+h_{m+1,m}\left[\begin{array}{c} 
            v_{m+1}^{(\mn)} \\
            v_{m+1}^{(p)}
        \end{array}\right]e_m^T
    \end{aligned}
\end{equation*}
using the definition of $\Tilde{A}$ as in \eqref{eq_4_assembled_matrix}. Assuming $\xi_m=\infty$, formula \cref{eq_3_3_approx} can be applied, which justifies the approximation
\begin{align*}
    e^{h_i\Tilde{A}}\left[\begin{array}{c} 
    c_0 \\
    e_p
\end{array}\right]\approx c_{\text{norm}}\left[\begin{array}{c} 
    V_m^{(\mn)} \\
    V_m^{(p)}
\end{array}\right]e^{h_iH_mK_m^{-1}}e_1\text{.}
\end{align*}
Equivalent to the derivation in \cref{s_4_1}, this finally leads to
\begin{align*}
    \sum_{k=0}^ph_i^k\varphi_k(-h_i\tA)\tc_k
    &\approx c_\text{norm}\left(\sum_{i=1}^m\tv_i\left[e^{h_iH_mK_m^{-1}}\right]_{i,1}\right)\text{.}
\end{align*}

In conclusion, it can be stated once again that no tensor matrix needs to be matricized and no tensor needs to be vectorized in the practical implementation of the method. This makes it possible to perform all relevant calculations in the \tttext\ format. The method derived here is called \rkexptt.

\section{Numerical Experiments}\label{c5}
In this section, the classical \cite{RK2EXPINT, KIOPS} and \tttext\ versions of \kiops\ and \rkexp\ are evaluated from various perspectives, particularly regarding runtime and accuracy. They are therefore used to solve the Allen-Cahn equation in $d\in\N$ dimensions \cite{RK2EXPINT, KIOPS, allen_cahn} and the heat equation.

\subsection{Problem setting}

Two equations are considered in the numerical experiments of this work. Let $\Delta$ denote the $d$-dimensional Laplace operator and $t\in[0,1]$. The first equation is the heat equation
\begin{equation*}
    \frac{\partial\tu(t)}{\partial t} = \Delta\tu(t)
\end{equation*}
which represents a simple model without any nonlinearities. Consequently, it can be solved directly by computing the matrix exponential, requiring no further $\varphi$-functions. The second equation under consideration is the Allen-Cahn equation
\begin{equation*}
    \frac{\partial\tu(t)}{\partial t}=\varepsilon^2\Delta\tu(t)+\tu(t)-(\tu(t))^3 = \varepsilon^2\Delta\tu(t)+\tu(t)(1-(\tu(t))^2)
\end{equation*}
where all multiplications are performed element-wise and $\varepsilon\in\R$. In the context of the \tttext\ format, the latter formulation of the Allen-Cahn equation is preferred because it reduces the number of \tttext\ truncations. The evaluation of the matrix algorithms will be based on the vectorized forms of both equations.

Depending on the geometry of the domain, the $d$-dimensional Laplace operator can be discretized in various ways. The finite difference discretization is presented in the following definition.

\begin{definition}[\cite{RK2EXPINT, lehrbuch_projector}]
    The real symmetric finite difference matrix $T_{n}$ of an equispaced triangulation of length $L$ is given by
    \begin{equation*}
        T_n=\frac{1}{h^2}\text{\normalfont{tridiag}\;\!}(-1,2,-1)\in\R^{n\times n}
    \end{equation*}
    where $h=\nicefrac{L}{n}$ denotes the spatial step size. The \fdm\ discretization of the $d$-dimensional Laplace operator $\Delta$ in matrix form is
    \begin{equation*}
        \Delta\approx A=\sum_{i=1}^dI_{n^{i-1}}\otimes_KT_n\otimes_KI_{n^{d-i}}\in\C^{n^d\times n^d}
    \end{equation*}
    where $\otimes_K$ denotes the Kronecker product, $I_k\in\C^{k\times k}$, $k\in\N$, is the identity matrix and $I_0:=1$. In the tensor case, the finite difference (\fdm) discretization of $\Delta$ reads
    \begin{equation*}
        \Delta\approx \tA=\sum_{i=1}^d\tI_{\mn_{(i-1)}}\otimes T_n\otimes\tI_{\mn_{(d-i)}}\in\C^{\mn\times\mn}
    \end{equation*}
    where $\otimes$ denotes the tensor product and $\tI_{\mn_{(k)}}\in\C^{\mn_{(k)}\times\mn_{(k)}}$ is the identity tensor for the mode set $\mn_{(k)}=(n,\ldots,n)\in\C^k$, $k\in\N$, $\mn:=\mn_{(d)}$ and $\tI_{\mn_{(0)}}:=1$.
\end{definition}

The \fdm\ discretization of the $d$-dimensional Laplace operator typically yields highly sparse matrices, providing a substantial computational advantage for sparse matrix methods. To evaluate performance on less sparse structures, the operator is also discretized on non-trivial geometric domains, specifically a cube, an annulus, and a flag, by utilizing methods from Isogeometric Analysis.

\subsection{Low-rank Isogeometric Analysis}
\label{s_5_2}

\newcommand*{\bigtimes}{\mathop{\raisebox{-.5ex}{\hbox{\huge{$\times$}}}}} 

Isogeometric Analysis (\iga), introduced in \cite{CAD}, bridges computer-aided design and finite element analysis by employing the same basis functions for both the parametrization of the computational domain and the approximation of the solution fields. Following this isoparametric concept, a geometry map $G \colon [0,1]^3 \to \Omega$ represents the computational domain $\Omega$ using a tensor-product B-spline or \nurbs\ space. To solve a partial differential equation, such as Poisson's equation, within a standard Galerkin framework, the discrete solution is searched in the same space or in a suitably refined version of it.

After pulling back the bilinear form from the physical domain $\Omega$ to the parametric domain $[0,1]^3$, the entries of the stiffness tensor $\mathbf{K} \in \mathbb{R}^{\mathbf{n}\times \mathbf{n}}$, with mode size ${\mn = (n_1, n_2, n_3)}$, are given by
\begin{equation*}
  (\mathbf{K})_{\mathbf{i},\mathbf{j}}
  =
  \sum_{k,l=1}^{3}
  \int_{[0,1]^3}
  \frac{\partial}{\partial \hat{x}_l} \beta_{\mathbf{i}}\left(\hat{x}\right)
  \frac{\partial}{\partial \hat{x}_k} \beta_{\mathbf{j}}\left(\hat{x}\right)
  q_{k,l}\left(\hat{x}\right)
  \,\mathrm{d}\hat{x},
\end{equation*}
where $\mathbf{i} = (i_1, i_2, i_3) \in \bigtimes^{3}_{d = 1} \left\{ 1, \ldots, n_d \right\}$ is a multi-index and $\beta_{\mathbf{i}} = \beta^{(1)}_{i_1} \beta^{(2)}_{i_2} \beta^{(3)}_{i_3}$ denotes the tensor-product B-spline basis function. The geometry-induced weight functions are defined by $\left(q_{k,l}\left(\hat{x}\right)\right)_{k,l=1}^{3} = \det\left(\nabla G\left(\hat{x}\right)\right) \left(\nabla G\left(\hat{x}\right)\right)^{-1} \left(\nabla G\left(\hat{x}\right)\right)^{-\top}$. Although the underlying B-spline basis functions possess a strict tensor-product structure, these weight functions are generally not separable, due to their rational dependence on the Jacobian of the geometry map. Consequently, the integrands of the stiffness matrix cannot, in general, be split into products of univariate factors. This leads to computationally expensive multivariate quadrature in the standard assembly process.

To overcome this bottleneck, Mantzaflaris et al.~\cite{angelos1} proposed approximating the non-separable weight functions by projecting them onto a sufficiently rich tensor-product B-spline interpolation space and compressing the resulting coefficient tensors using low-rank formats, such as the higher-order \svd\ (\hosvd). This yields an approximation of the form
\begin{equation*}
    \text{\small$
  q_{k,l}\left(\hat{x}\right)
  \approx
  \left\langle \tilde{\mathbf{V}}_{k,l},\tilde{\mathbf{B}}\left(\hat{x}\right) \right\rangle_F
  = \left\langle
    \sum_{r=1}^{R_{k,l}}
  \bigotimes^3_{d = 1} 
  v_{k,l,r}^{(d)}
  , \bigotimes^3_{d = 1} \tilde{B}^{(d)} \left(\hat{x}_d \right) \right\rangle_F
  =
  \sum_{r=1}^{R_{k,l}}
  \prod_{d=1}^{3}
  \left(
    v_{k,l,r}^{(d)}
    \cdot
    \tilde{B}^{(d)}\left(\hat{x}_d\right)
  \right),
  $}
\end{equation*}
where $v_{k,l,r}^{(d)}\in\mathbb{R}^{\tilde{n}_d}$ are the factors of the summands in a canonical polyadic notation and $\tilde{B}^{(d)}\left(\hat{x}_d\right) \in \mathbb{R}^{\tilde{n}_d}$ contain the univariate basis function in mode $d = 1,2,3$ of the projection B-spline space. Once the weight functions are replaced by their low-rank separable interpolants, the three-dimensional integrals decouple into products of univariate integrals. Consequently, the global stiffness matrix can be efficiently approximated as a sum of Kronecker products of small univariate matrices,
\begin{equation*}
  \mathbf{K}
  \approx
  \sum_{k,l=1}^{3}
  \sum_{r=1}^{R_{k,l}}
  \bigotimes_{d=1}^{3} K_{k,l,r}^{(d)},
\end{equation*}
where
\begin{gather*}
  K_{k,l,r}^{(d)}
  =
  \int_{0}^{1}
    \left( 
  \delta(k,d) B^{(d)}\left(\hat{x}_d\right)
  \otimes
  \delta(l,d) B^{(d)}\left(\hat{x}_d\right)
  \right)
   \left(
    v_{k,l,r}^{(d)}
    \cdot
    \tilde{B}^{(d)}\left(\hat{x}_d\right)
  \right)
  \,\mathrm{d}\hat{x}_d, \\
  \delta(k,d) f
  =
  \begin{cases}
    \dfrac{\partial f}{\partial \hat{x}_d}, & k=d,\\[0.6em]
    f, & k\neq d.
  \end{cases}
\end{gather*}
Interpolating the weight functions in this way typically yields full coefficient tensors. B\"unger et al.~\cite{BuengerDolgovStoll:2020} extended this methodology by directly computing the coefficient tensors in the TT format. By exploiting the natural tensor-product structure of both the interpolating B-spline basis and the corresponding grid of Greville points, the interpolation problem can be formulated as a linear system that is solved using the Alternating Minimal Energy (\amen) solver \cite{amen}. This avoids the explicit construction of the dense interpolation tensors $\mathbf{V}_{k,l}$, keeps the TT ranks moderate, and assembles the stiffness matrix efficient and scalable. This approach was further extended in \cite{Riemer2025} to multi-patch geometries and in \cite{Riemer2026} to locally refined solution spaces based on (truncated) hierarchical B-splines. 

In this work, the method presented in \cite{BuengerDolgovStoll:2020} is used to assemble the considered stiffness tensors low-rank.

\subsection{Experiments}

\begin{figure}[htbp]
\centering
\begin{subfigure}[t]{0.29\textwidth}
\centering
\includegraphics[width=\linewidth]{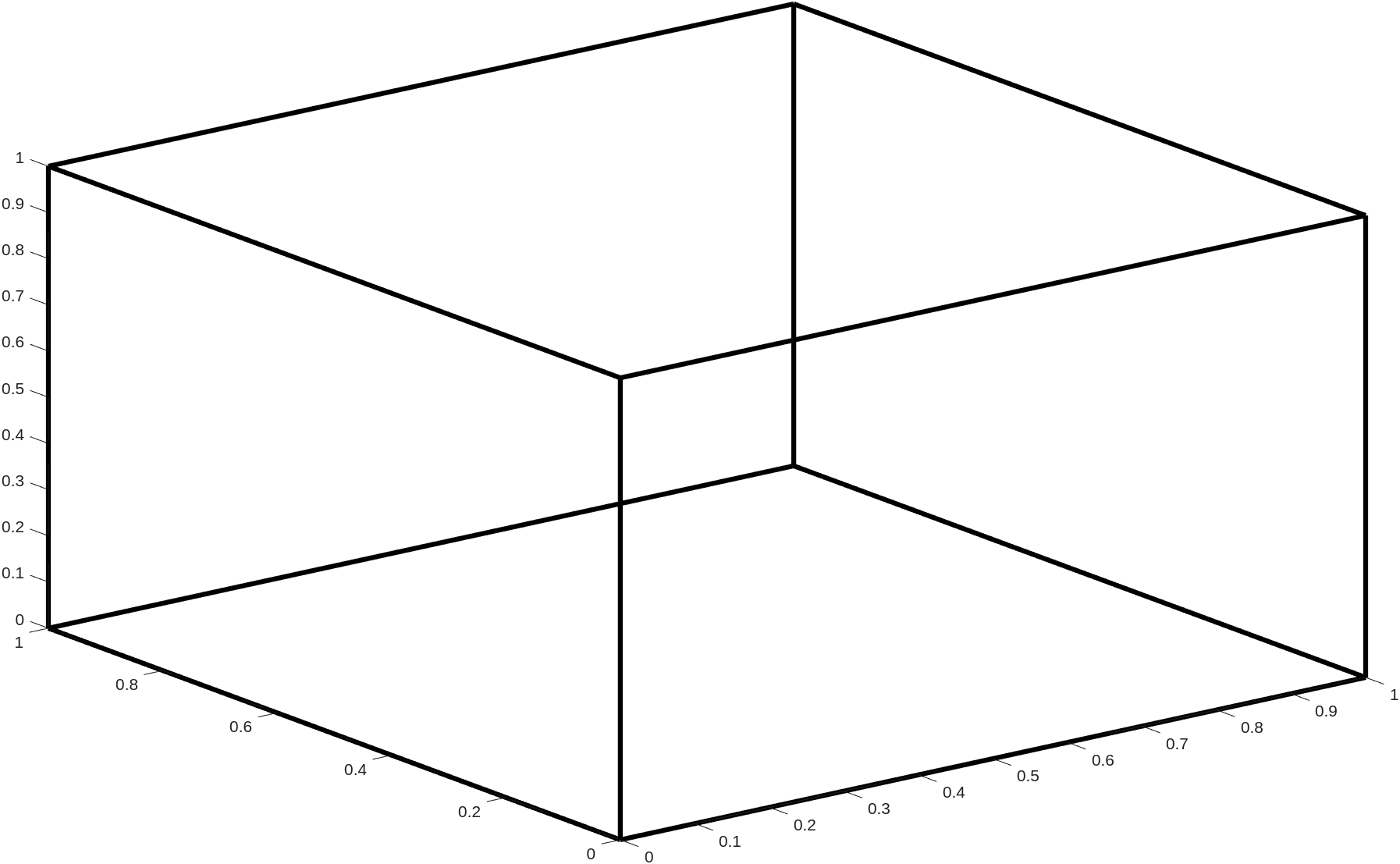}
\label{fig:1_a}
\end{subfigure}\hfill
\begin{subfigure}[t]{0.29\textwidth}
\centering
\includegraphics[width=\linewidth]{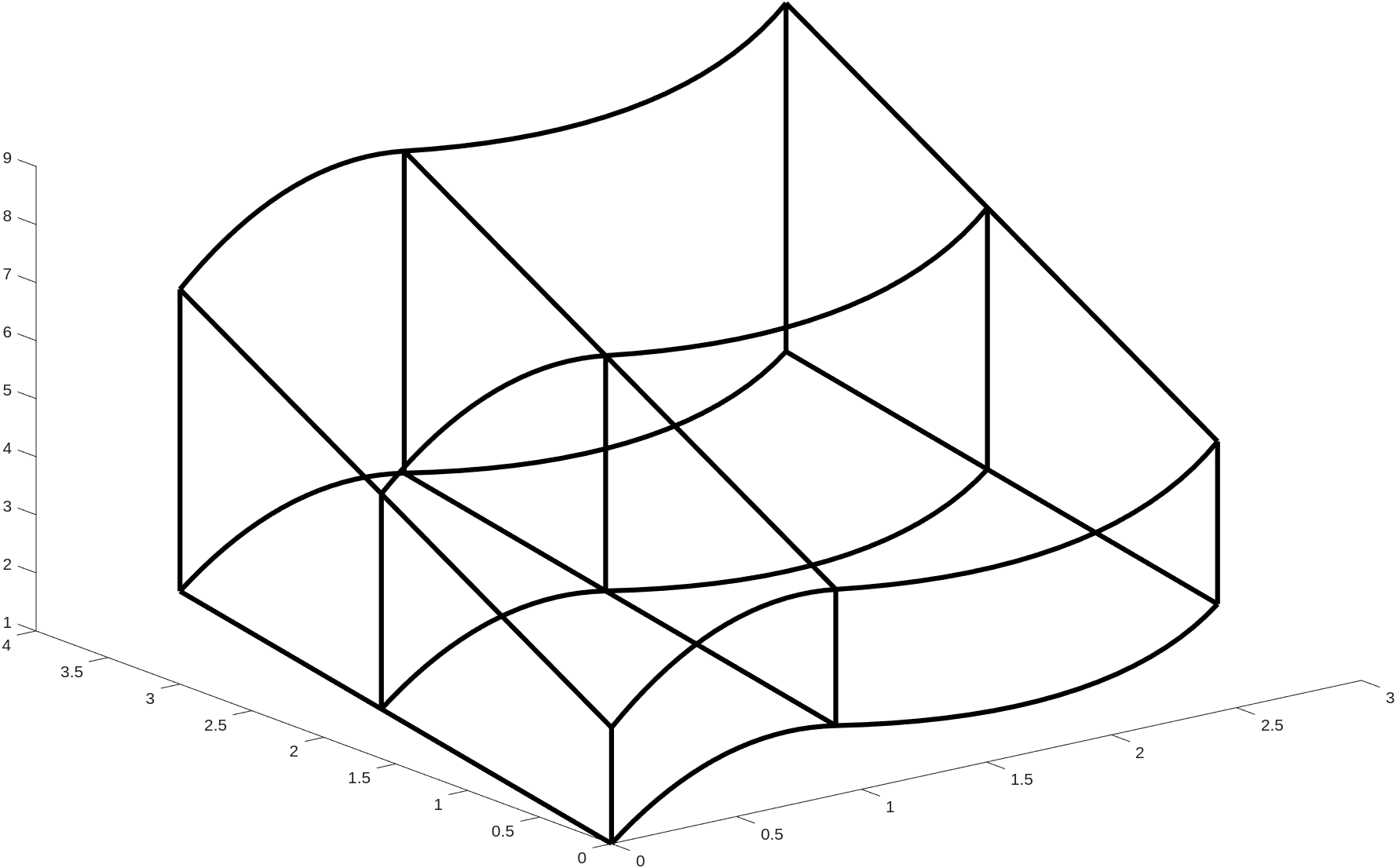}
\label{fig:1_b}
\end{subfigure}\hfill
\begin{subfigure}[t]{0.29\textwidth}
\centering
\includegraphics[width=\linewidth]{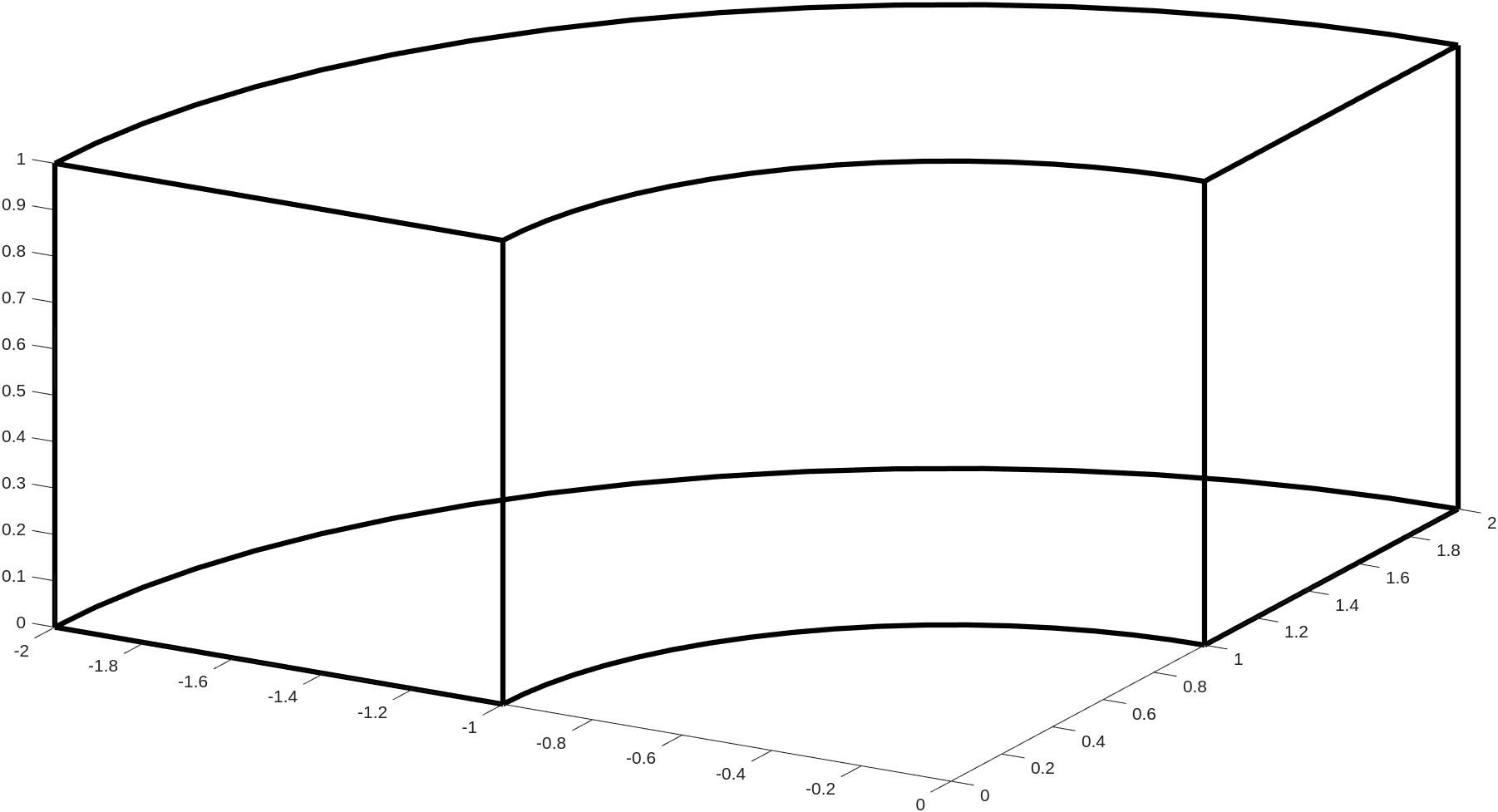}
\label{fig:1_c}
\end{subfigure}
\\[5pt]\caption{The two B-spline geometries of a cube (left) and a \emph{Thick} flag (center) and the \nurbs\ geometry of a quarter annulus (right) used in the numerical experiments.\\[-25pt]}
\label{fig:1}
\end{figure}

The default setting of the experiments is as follows: The number of dimensions is $d=3$. For the \fdm\ discretization of $\Delta$ as the Laplace operator, $n=100$ is the default size of each dimension, i.e. the number of grid points in the three-dimensional case is $n^3$. For the \iga\ setting, experiments are conducted on two B-spline geometries and one \nurbs\ geometry, which are depicted in \cref{fig:1}. After interpolating the weight functions according to the procedure in \cref{s_5_2}, the stiffness tensor of the cube and the quarter annulus is assembled. This assembly utilizes a multivariate tensor-product B-spline basis with $n_1=n_2=n_3=30$ basis functions. For the \emph{thick} flag configuration, the basis functions are increased to ${n_1=n_2=60}$ and $n_3=30$. By default, B-splines of degree 5 are employed.

The time interval $T=[0,1]$ is divided into 20 equally sized subintervals with step size $h=0.05$. The methods \kiops, \kiopstt, \rkexp\ and \rkexptt\ are given an accuracy of $\text{tol}=10^{-8}$ with which the solution is to be calculated. A reference solution is calculated by \matlab s ode15s solver with an accuracy of $10^{-10}$. For \tttext\ operations like rounding, the tolerance is empirically set to $\text{tol}_{\text{\tttext}}=\text{tol}$. The heat equation is solved by evaluating the action of the exponential $e^{\varepsilon h\tA}$ on $\tu_i$, i.e. $\tu_{i+1}=e^{\varepsilon h\tA}\tu_i$, directly, since this already yields the exact solution. To solve the Allen-Cahn equation, the exponential integrator \krogstad\ is used, as it demonstrated the best trade-off between accuracy and runtime in preliminary simulations. 
Based on further test runs, the default pole set for the rational methods is chosen to consist of 30 complex poles optimized by the \rkfit\ algorithm \cite{rkfit} for experiments using the \fdm\ discretization of $\Delta$ and of a single pole $\xi=31400$ with a multiplicity of 72 for problems using the \iga\ discretization. The pole sets are taken from \rkexp\ \cite{RK2EXPINT, rk2expint_repo}.
For solving their respective linear systems, \rkexp\ employs the \agmg\ algorithm \cite{RK2EXPINT, agmg} and \rkexptt\ utilizes the \amen\ solver \cite{amen}.

The source code for reproducing the results and figures of this section is accessible at \url{https://github.com/riweig/ttexpint}.

\subsubsection*{Pole selection} 

The same pole sets as in \cite{RK2EXPINT} are considered here. These are three real pole sets consisting of one, two, or four distinct repeated poles, and two complex pole sets. More precisely, one of the complex pole sets is derived from the rational best approximants to the function $e^{z}$ on the negative real semi-axis, using the roots of the denominator polynomial of the approximation as poles, while the other consists of complex-valued poles selected by the \rkfit\ algorithm \cite{guettel_toolbox, rkfit} subject to the condition of having positive real parts. The real pole sets comprise 72 elements each, while the complex versions feature 30 poles each.

\Cref{s_4_2} notes that early \rkexptt\ implementations suffered from stability issues caused by the \amen\ solver when handling complex poles. To resolve these issues, the final algorithm treats the real and imaginary parts of the TT tensors separately, thereby executing all operations in the real domain. Experimental results show that this procedure is successful and provides increased computational stability. Particularly, the accuracy of the final version of \rkexptt\ is at the same level as \rkexp, even when complex poles are considered.

\begin{figure}[htbp]
    \centering
    \begin{tabular}{c@{\hspace{3pt}}c}
        \includegraphics[width=4.1cm]{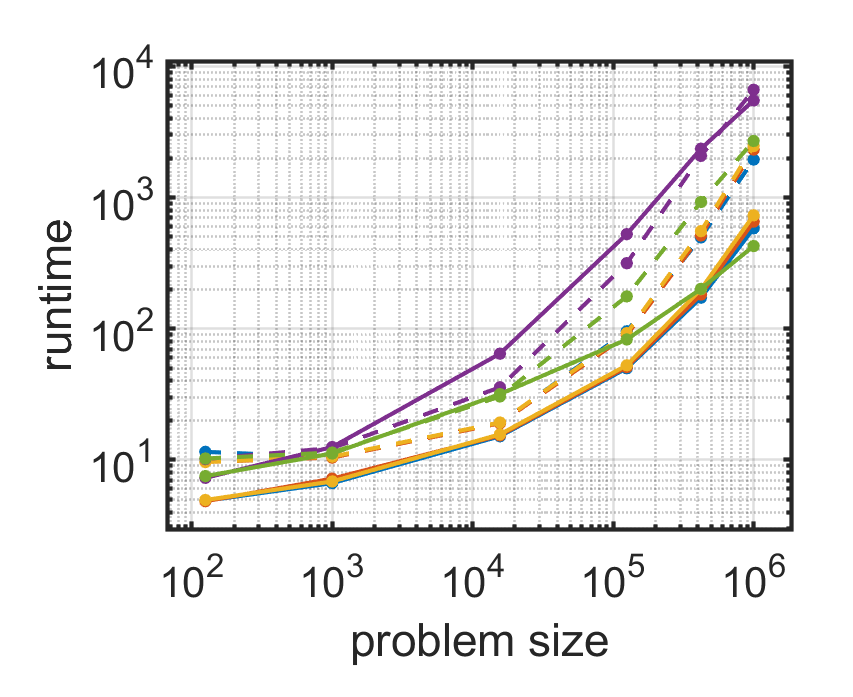} &
        \includegraphics[width=4.1cm]{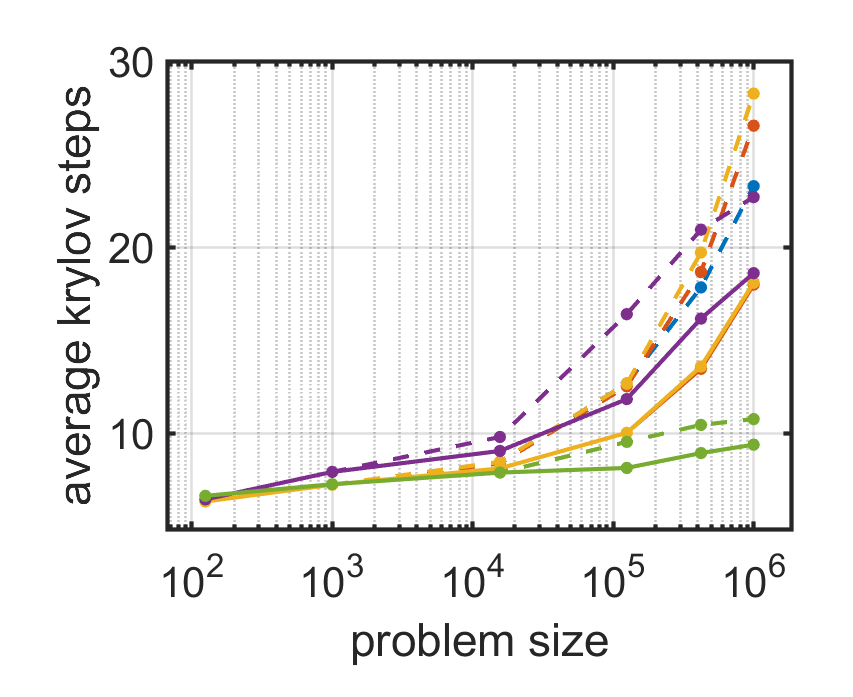} \\
        \multicolumn{2}{c}{\includegraphics[height = .5cm]{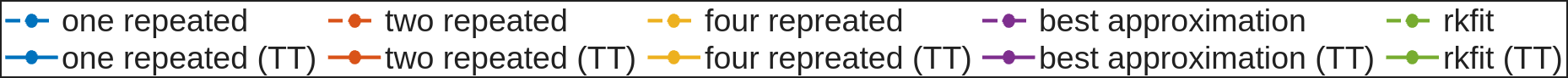}} \\[-9pt] \\ \\
    \end{tabular}
    \caption{Comparison of the pole sets on the algorithms \rkexp\ and \rkexptt. The sets are evaluated on the three-dimensional Allen-Cahn equation on a \fdm\ grid with a varying number of grid points.\\[-25pt]}
    \label{fig:pol1}
\end{figure}

This is particularly evident when considering problems involving the discretization of the equation on a \fdm\ grid. Indeed, a glance at \cref{fig:pol1} confirms the findings of~\cite{RK2EXPINT}: while the repeated real pole sets exhibit good runtimes, the results for the pole set of the best rational approximations are rather poor. Regarding the average number of Krylov steps, the poles optimized by the \rkfit\ algorithm yield the lowest values, also showing a significantly smaller increase than all other sets. As seen with the TT version, this also impacts the runtime for the largest problem sizes. These results justify the default use of \rkfit\ poles for problems discretized on \fdm\ grids.

However, for problems involving \iga\ discretization, the average Krylov iteration numbers for all pole sets remain constant at the level of the specified minimum. In some cases the iteration counts are even below, due to early breakdowns upon finding an optimal basis. For this reason, the increased computational complexity associated with complex poles leads to longer runtimes. This explains the default choice of a repeated real pole when considering problems with \iga\ discretization.

\subsubsection*{Accuracy}

\begin{figure}[htbp]
    \centering
    \begin{tabular}{c}
        \includegraphics[width=4.1cm]{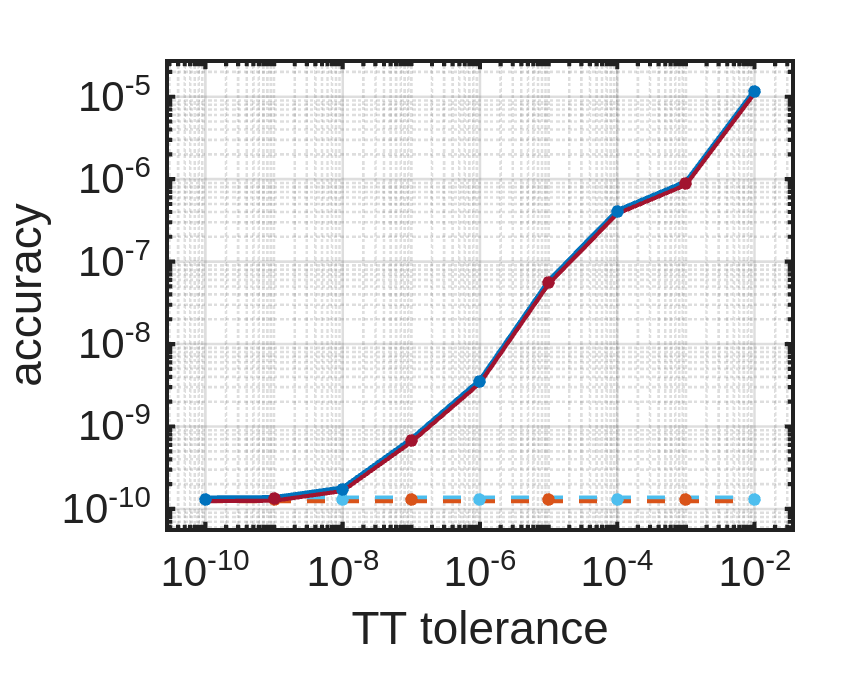} \\
        \includegraphics[height = .26cm]{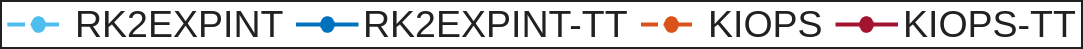} \\[-9pt] \\ \\
    \end{tabular}
    \caption{Comparison of the algorithms for evaluating the matrix exponential. The methods are evaluated on the Allen-Cahn equation on an \iga\ flag. The solution is computed with different \tttext\ operation precision.\\[-25pt]}
    \label{fig:acc1}
\end{figure}

The accuracy is measured by the infinity norm of the difference between the reference solution and the solution computed via \kiops, \kiopstt, \rkexp, or \rkexptt. Assuming that inaccuracies in the solutions accumulate as iterations progress, the computed solution is considered only at the last point of the time interval. In all experiments, the accuracy of the \tttext\ methods is on par with that of the matrix algorithms, being sometimes slightly better and sometimes slightly worse. Notably, the deviation from the reference solution does not exceed $10\text{tol}$, where $\text{tol}$ denotes the tolerance used for the computation of the matrix exponential. As approximation inaccuracies for the non-linear part of the problem equation are expected in exponential integrators, a factor of 10 is reasonable. Only when step sizes larger than the default are considered or the \tttext\ truncation is applied with lower precision does the error exceed $10\text{tol}$. This is expected, as the step size directly affects the approximation quality, especially regarding the non-linear component. The dependence of the solution quality of the \tttext\ methods on the \tttext\ operation tolerance $\text{tol}_{\text{\tttext}}$ is shown in \cref{fig:acc1}. It also illustrates that $\text{tol}_{\text{\tttext}}$ should not be chosen too generously.

\subsubsection*{Runtime}

\begin{figure}[htbp]
    \centering
    \begin{tabular}{c@{\hspace{3pt}}c@{\hspace{3pt}}c}
        \includegraphics[width=4.1cm]{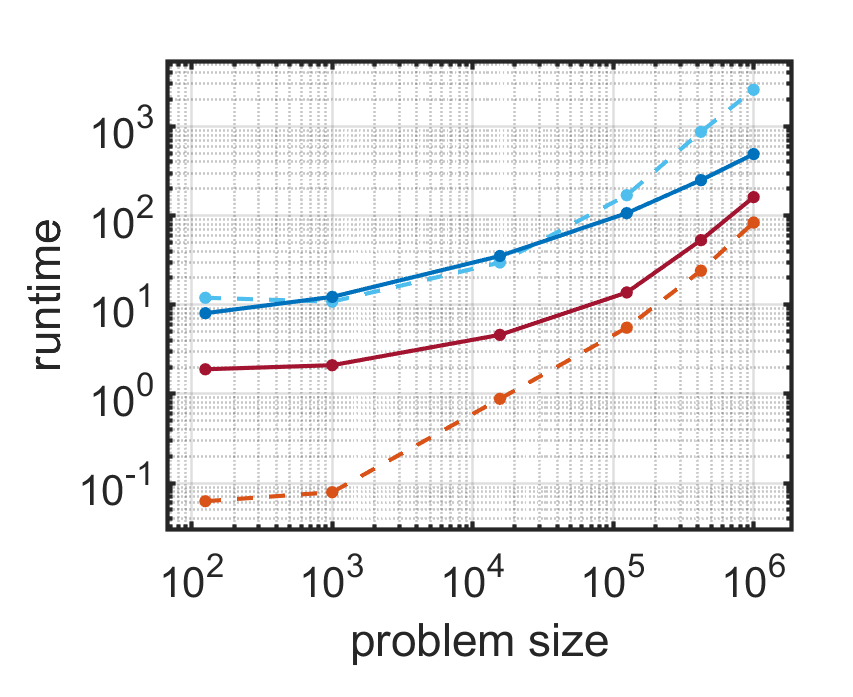} &
        \includegraphics[width=4.1cm]{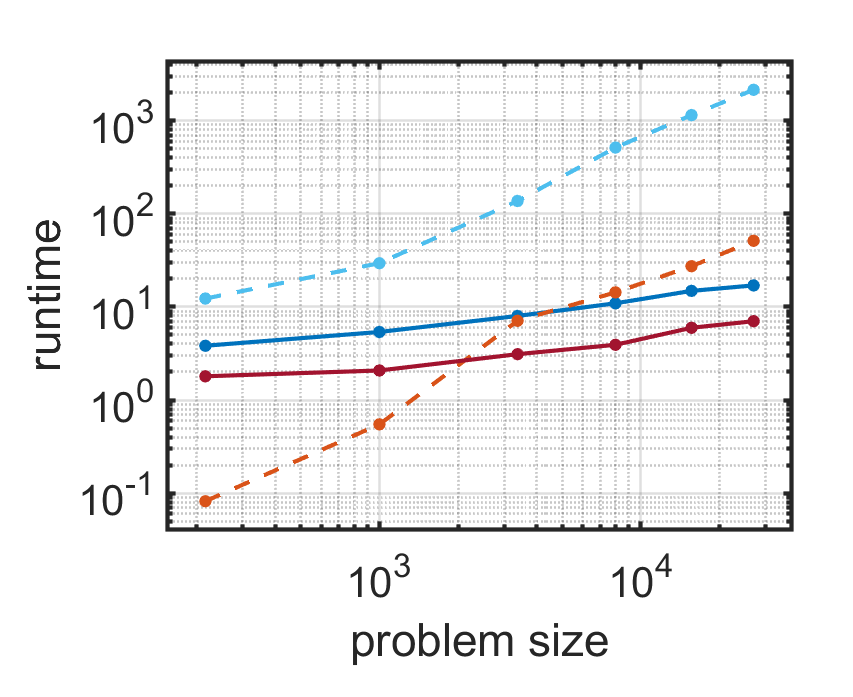} &
        \includegraphics[width=4.1cm]{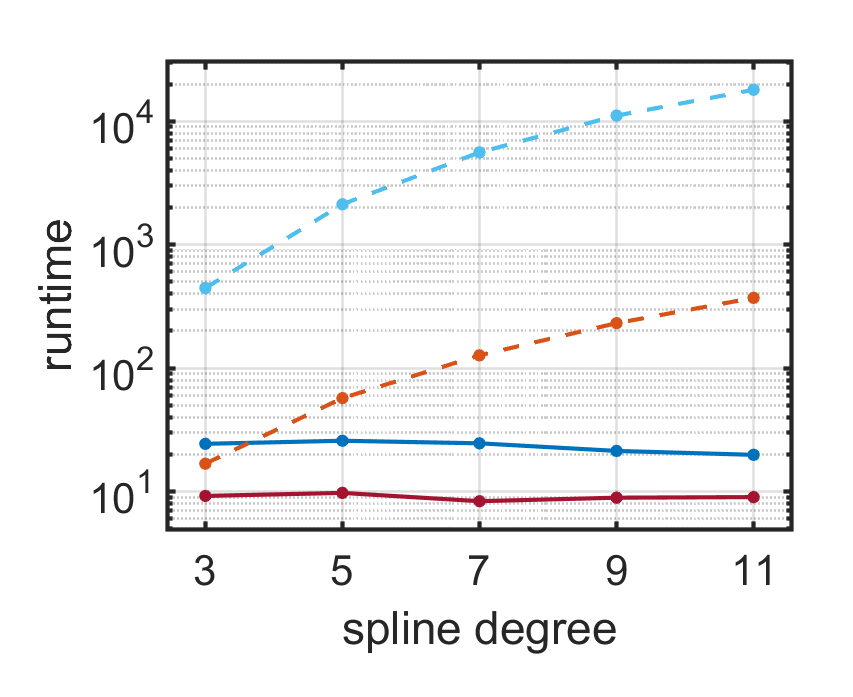} \\
        \multicolumn{3}{c}{\includegraphics[height = .26cm]{figures/experiments/leg_normal.png}} \\[-9pt] \\ \\
    \end{tabular}
    \caption{Comparison of the algorithms for evaluating the matrix exponential. The methods are evaluated on the three-dimensional Allen-Cahn equation on a \fdm\ grid with a varying number of grid points (left), on an \iga\ cube with varying numbers of basis functions (center) and on an \iga\ annulus with varying spline degree (right).\\[-25pt]}
    \label{fig:time1}
\end{figure}

Inference time is a critical criterion in practical use. The first two experiments in \cref{fig:time1} show that the runtime curves of the \tttext\ methods appear to at least approach that of \kiops\ with increasing problem size, and even outperform it depending on the problem. However, for the \fdm\ discretization of the Allen-Cahn equation, evaluating larger problems is required to prove that the \tttext\ version surpass \kiops. This was unfortunately unfeasible due to resource constraints. In contrast, for the IgA discretizations, it quickly becomes apparent that the \tttext\ methods outperform the matrix algorithms in terms of runtime as the problem size increases.

Since the Laplace operator matrix is highly sparse when using the \fdm\ discretization, \kiops\ and \rkexp\ seem to have a substantial advantage, making them hard to surpass in this scenario. To test this hypothesis, the third experiment in \cref{fig:time1} successively increases the spline degree of the IgA discretization, leading to denser matrices. As anticipated, this drives up the runtimes of the matrix algorithms. Meanwhile, the execution time of the \tttext\ algorithms remains almost untouched. The reason is that the \tttext\ ranks of the discretized linear operators stay nearly constant as the spline degree grows.

Notably, the average Krylov iteration numbers of \kiops\ and \kiopstt, as well as of \rkexp\ and \rkexptt, show very similar trends. This behavior was expected based on theoretical foundations, though regular \tttext\ roundings introduced an element of uncertainty. Because the average number of Krylov steps evolve similarly, the runtime benefits of the \tttext\ methods can be directly credited to more efficient computations using \tttext\ objects. Consequently, varying convergence speeds were ruled out as a contributing factor.

\begin{figure}[htbp]
    \centering
    \begin{tabular}{c@{\hspace{3pt}}c@{\hspace{3pt}}c}
        \includegraphics[width=4.1cm]{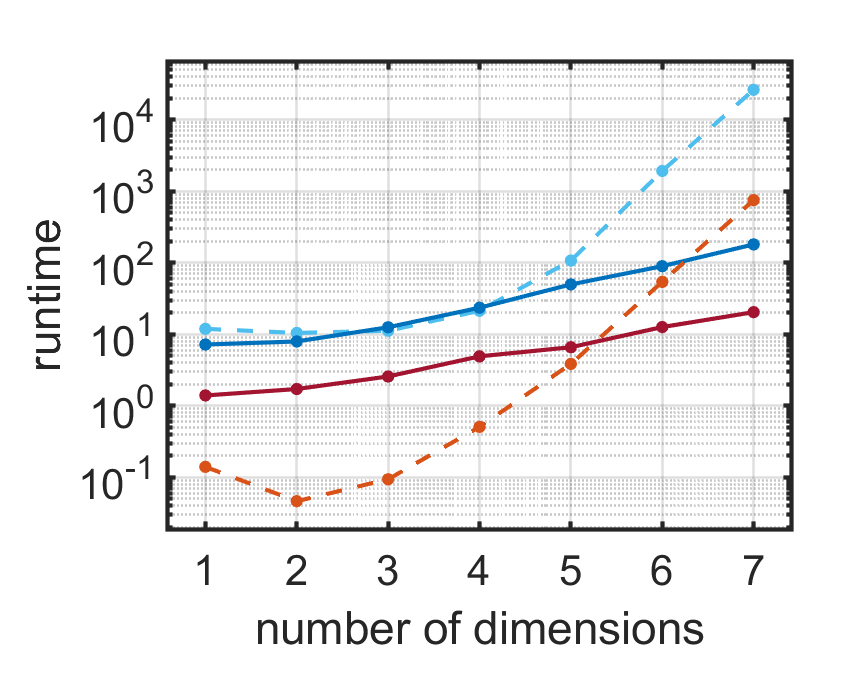} &
        \includegraphics[width=4.1cm]{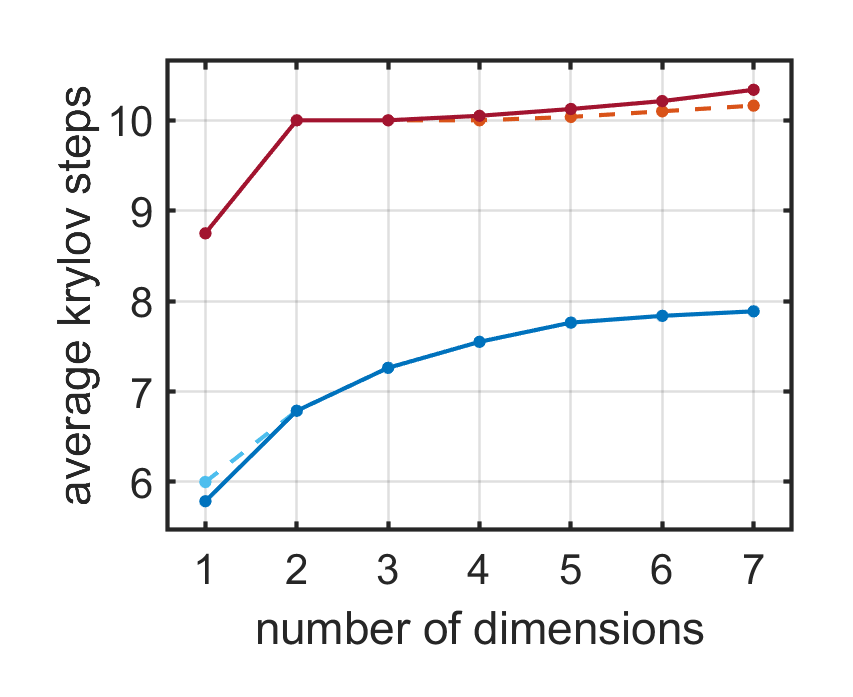} &
        \includegraphics[width=4.1cm]{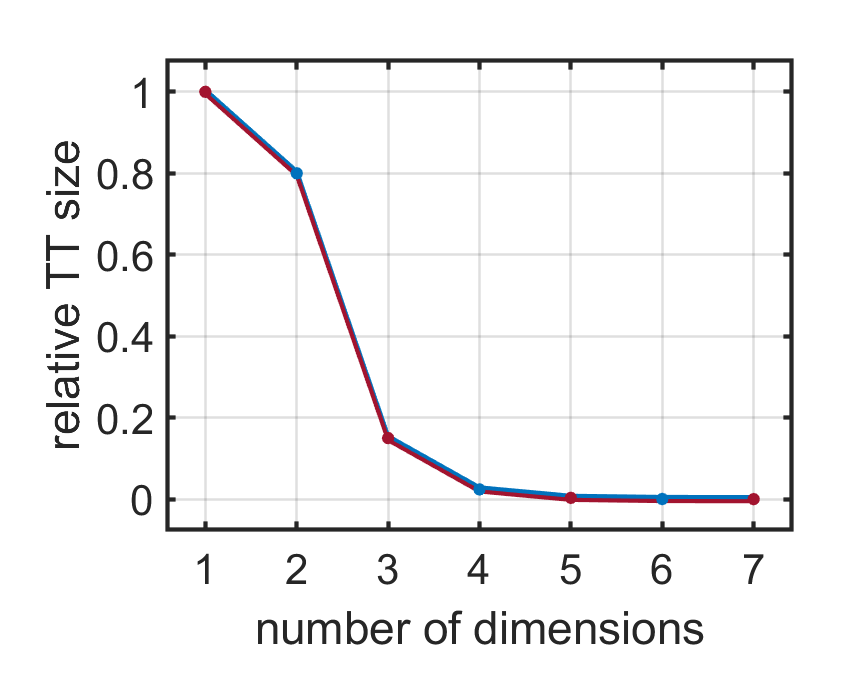} \\
        \multicolumn{3}{c}{\includegraphics[height = .26cm]{figures/experiments/leg_normal.png}} \\[-9pt] \\ \\
    \end{tabular}
    \caption{Comparison of the algorithms for evaluating the matrix exponential. The methods are evaluated on the Allen-Cahn equation on a \fdm\ grid with increasing numbers of dimensions. The choice of $n=10$ (instead of the default choice of $n=100$) for each dimension implies a total of $10^{d}$ grid points.\\[-25pt]}
    \label{fig:time2}
\end{figure}

The experiment in \cref{fig:time2} scales the number of dimensions to increase the problem size, instead of increasing the number of grid points or B-Spline basis functions as done elsewhere. This specifically targets the curse of dimensionality, which tensor decompositions aim to break. As expected, the matrix algorithms quickly show clear exponential runtime growth. Conversely, the growth for \tttext\ methods is much less steep, enabling them to significantly outperform classical approaches at higher dimensionalities. The iteration numbers again develop very similarly, although it should be noted that for \kiops\ and \kiopstt\ they are close to the specified minimum. While it is possible to increase the number of dimensions even further for the \tttext\ algorithms, the matrix methods hit the memory limits of the infrastructure used. A look at the relative \tttext\ size of the respective final solution, which is defined as the number of elements stored by a \tttext\ tensor relative to the number of entries the tensor would contain in full format, also reveals that using a tensor decomposition is worthwhile for high-dimensional problems.

\begin{figure}[htbp]
    \centering
    \begin{tabular}{c@{\hspace{3pt}}c}
        \includegraphics[width=4.1cm]{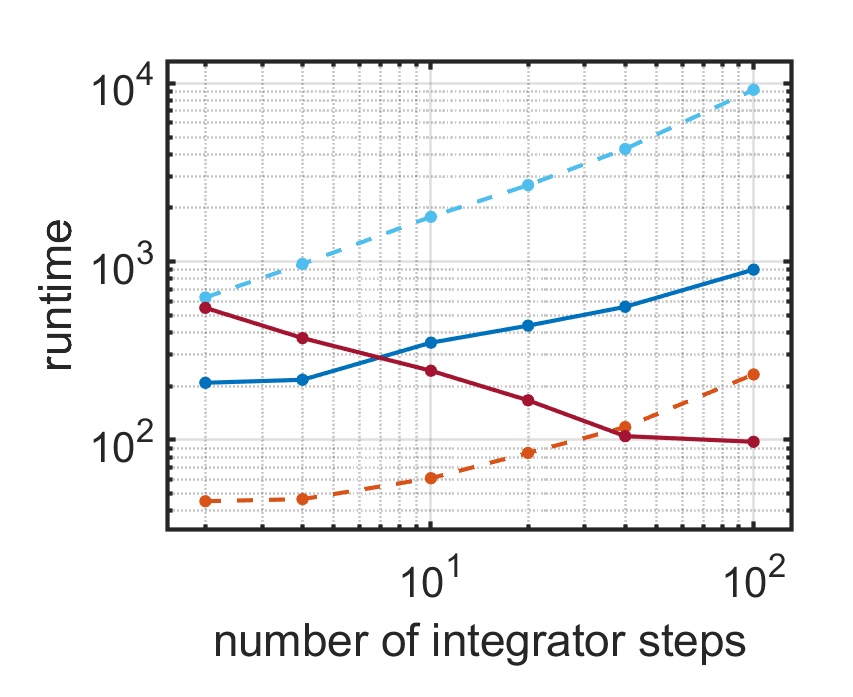} &
        \includegraphics[width=4.1cm]{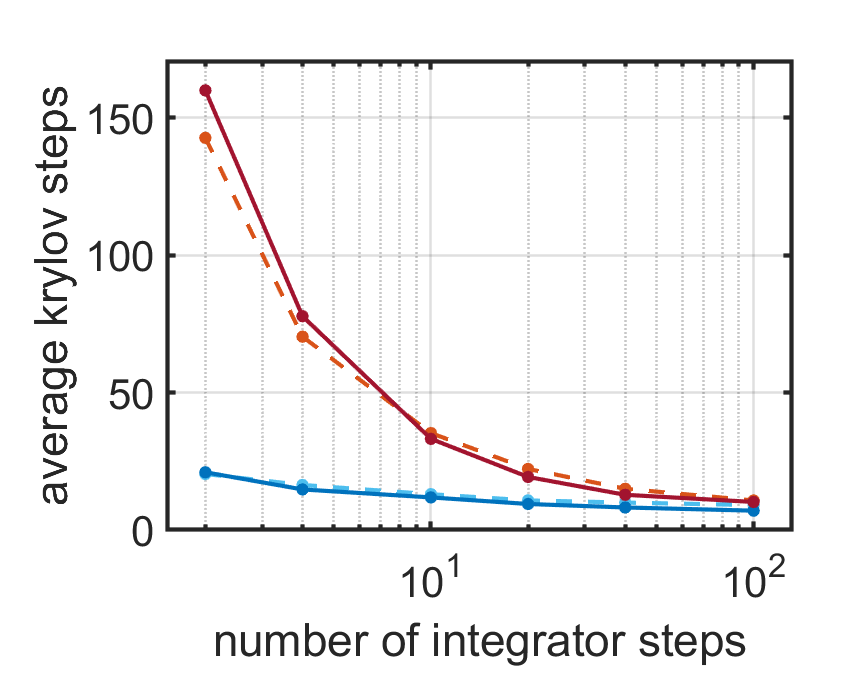} \\
        \multicolumn{2}{c}{\includegraphics[height = .26cm]{figures/experiments/leg_normal.png}} \\[-9pt] \\ \\
    \end{tabular}
    \caption{Comparison of the algorithms for evaluating the matrix exponential. The methods are evaluated on the three-dimensional Allen-Cahn equation on a \fdm\ grid. The step size of the exponential integrator, and therefore the total number of time steps, is varied in this experiment.\\[-25pt]}
    \label{fig:time3}
\end{figure}

An interesting phenomenon is observed for the three-dimensional Allen-Cahn equation with \fdm\ discretization when varying the number of time steps. While the inference times of the algorithms \kiops, \rkexp, and \rkexptt\ increase as expected when refining the time interval, the inference time decreases for \kiopstt. Actually, the runtime was expected to increase as the step size decreases due to the growing number of evaluations of the exponential function. The phenomenon is visualized in \cref{fig:time3}.

\begin{figure}[htbp]
    \centering
    \begin{tabular}{c@{\hspace{3pt}}c}
        \includegraphics[width=4.1cm]{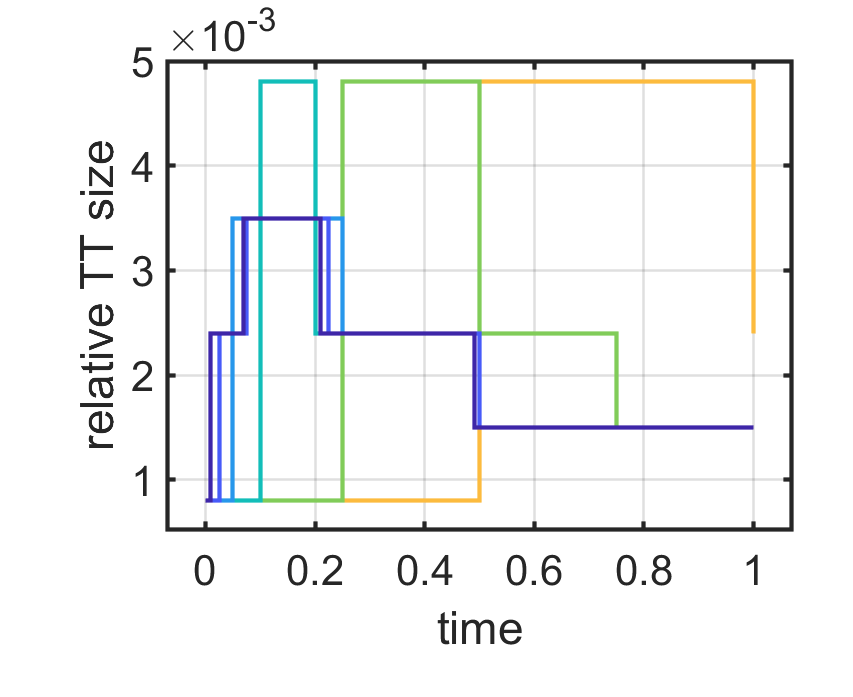} &
        \includegraphics[width=4.1cm]{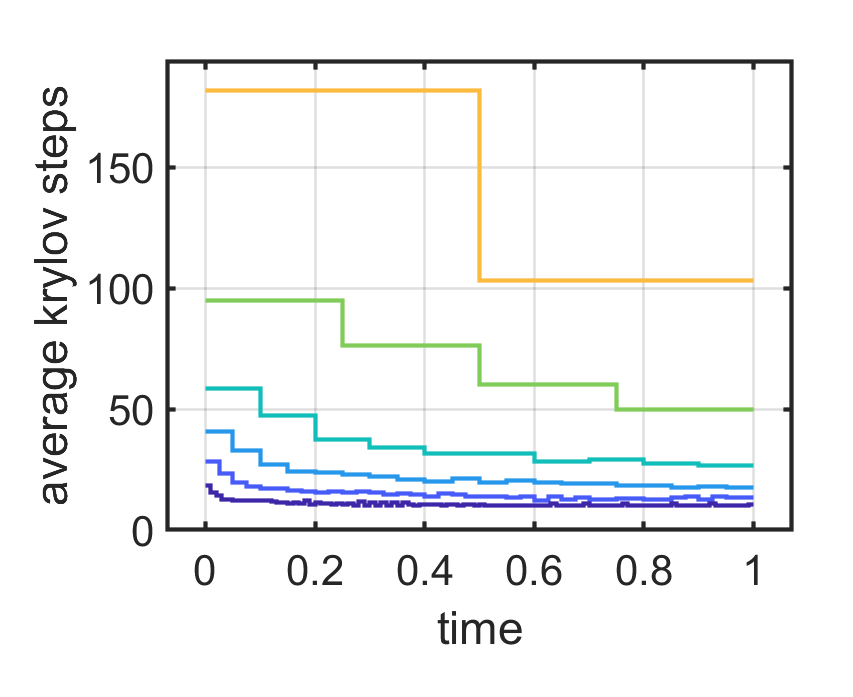} \\
        \multicolumn{2}{c}{\includegraphics[height = .26cm]{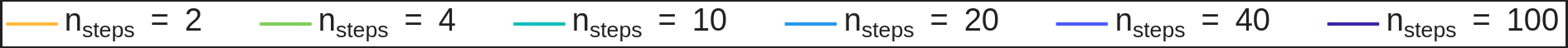}} \\[-9pt] \\ \\
    \end{tabular}
    \caption{Comparison of the algorithms for evaluating the matrix exponential. The methods are evaluated on the three-dimensional Allen-Cahn equation on a \fdm\ grid. The step size of the exponential integrator, and therefore the total number of time steps $n_{\text{steps}}$, is varied in this experiment.\\[-25pt]}
    \label{fig:time4}
\end{figure}

The explanation for this phenomenon is provided by a combination of various factors. First, it should be noted that the rapidly decreasing average number of Krylov iterations cannot be the sole justification, since the phenomenon of a decreasing runtime does not occur for \kiops\ in the matrix case. Furthermore, the rank evolution cannot be primarily responsible for the decreasing runtime either, because, as can be seen in \cref{fig:time4}, the ranks develop very similarly for 200 or more time steps. The key is provided by examining the evolution of the average Krylov iterations per time step in \cref{fig:time4}. These tend to decrease, and the different curves do not appear to intersect, or do so only marginally. Moreover, the number of Krylov steps, categorized by the number of time steps, differs much more significantly at the beginning than at the end. Considering the fact that the \tttext\ ranks tend to be high at the beginning of the time interval and low at the end, the phenomenon of decreasing runtimes for finer time intervals can be explained. However, this is presumably an exceptional case rather than the rule.

\subsubsection*{\tttext\ rank evolution}

The evaluation of the experiments shows no surprises regarding the rank evolution. The ranks evolve differently depending on the considered equation or discretization. For instance, after an initial increase, the ranks tend to decrease for the Allen-Cahn equation with \fdm\ discretization, while they increase or remain fairly constant in other scenarios. In the first time step, a significant increase in the ranks is always observed, which does not occur to this extent in subsequent time steps. This suggests that the rank of the initial conditions is artificially low. As expected, the relative \tttext\ ranks of the solutions decrease with increasing problem size, decreasing spline degree, and decreasing step size. Furthermore, the ranks for \kiopstt\ and \rkexptt\ evolve very similarly, with only minor deviations from one another.

Unfortunately, no clear indication could be found as to whether low ranks of the linear operator matter more than low ranks of the initial condition and solutions. However, because the tensor matrix has double the number of modes compared to the initial condition and solution tensors, it is likely that low \tttext\ ranks of the linear operator are far more important. Nevertheless, for the initial and solution tensors, there is also a dependency regarding the non-linear part of the initial value problem.

\section{Conclusions}\label{c6}
This work combined the advantages of the Tensor Train (\tttext) format \cite{tt_oseledets, TT} for storing and processing large amounts of data with solving stiff differential equations. The main contribution is the further development of \kiops\ \cite{KIOPS} and \rkexp\ \cite{RK2EXPINT} for use with the \tttext\ format, from which the methods \kiopstt\ and \rkexptt\ were derived. Furthermore, it was shown that the theoretical foundations of the matrix versions can also be transferred to the \tttext\ methods, thus ensuring the validity of convergence and accuracy properties. Additionally, some explicit exponential Runge-Kutta integrators \cite{RK2EXPINT, etd3rk, krogstad4, sw2} were adapted for use with the \tttext\ decomposition. However, it is emphasized that the use of the developed methods for evaluating tensor matrix exponential functions is not limited to explicit exponential Runge-Kutta integrators.

The theoretical results were verified by numerical experiments, which confirmed two important properties of \kiopstt\ and \rkexptt: accuracy at the same level as the matrix versions and runtime advantages. As expected, the runtime advantages were scenario-dependent. While \kiops\ and \rkexp\ continued to demonstrate their strength in very sparse stiffness matrices, they were outperformed by their \tttext\ counterparts in low-rank scenarios and with very large problems. In particular, \kiopstt\ performed often better than all other methods in these scenarios and is therefore recommended for use.

However, there is still potential for optimization, especially regarding the rational Krylov methods used. In particular, the question of an optimal pole selection strategy remains the subject of further research. One possibility would be to extend the \rkfit\ method \cite{rkfit}, which was developed to determine optimal poles in the matrix case, to tensors in the \tttext\ format.

The \tttext\ rounding procedure was implemented rather arbitrarily after up to five sums. However, a method that dynamically decides whether to round in the next step based on rank growth is also conceivable. Especially the automatic selection of an optimal \tttext\ tolerance is seen as a further optimization approach, possibly even an adaptive method for automatically adjusting the tolerance as needed.





\bibliographystyle{siamplain}
\bibliography{references}
\end{document}